\documentclass[preprint,12pt]{elsarticle}

\usepackage{mathrsfs}
\usepackage{amsfonts}
\usepackage{dsfont}
\usepackage{multirow}
\usepackage{bbold}
\usepackage{amssymb}
\usepackage{float}
\usepackage{placeins}
\usepackage[utf8]{inputenc}
\usepackage[margin=2cm]{geometry}
\usepackage{hyperref}

\usepackage{tikz}
\usetikzlibrary{backgrounds,arrows, shapes,decorations.markings,positioning}
\usetikzlibrary{plotmarks,calc,fadings,decorations.pathreplacing,decorations.pathmorphing}
\tikzset{%
  >=latex,
  inner sep=0pt,%
  outer sep=2pt,%
  mark coordinate/.style={inner sep=0pt,outer sep=0pt,minimum size=3pt,
    fill=black,circle}%
}

\tikzset{-dot-/.style={decoration={
  markings,
  mark=at position #1 with {\fill circle (2.0pt);}},postaction={decorate}}} %%% in this line added a ;

\usepackage{pgfplots}
\pgfplotsset{compat=1.13}
\usepackage{amsmath}
\usepackage{caption}
\usepackage{subcaption}
\usepackage{cleveref}
\newcommand{\R}{\mathbb{R}}

\newcommand{\diff}[1]{{\mathrm{d}{#1}}}

\newcommand{\polygon}[2]{%
  let \n{len} = {2*#2*tan(360/(2*#1))} in
 ++(0,-#2) ++(\n{len}/2,0) \foreach \x in {1,...,#1} { -- ++(\x*360/#1:\n{len})}}

\definecolor{darkgreen}{rgb}{0.0, 0.5, 0.0}
\DeclareMathOperator*{\argmin}{arg\,min}
\newtheorem{theorem}{Theorem}

\newtheorem{remark}[theorem]{Remark}

\journal{Elsevier}

\begin{document}

\begin{frontmatter}



\title{Very high order treatment of embedded curved boundaries in compressible flows: ADER discontinuous Galerkin with a space-time Reconstruction for Off-site data }


\author[1]{Mirco Ciallella\corref{cor1}}
\author[2]{Stephane Clain}
\author[3]{Elena Gaburro}
\author[3]{Mario Ricchiuto}

\address[1]{\'Ecole Nationale Sup\'erieure d'Arts et M\'etiers, Institut de M\'ecanique et d'Ing\'enierie, 33400 Talence, France}
\address[2]{CMUC –Mathematical centre, Coimbra University, Coimbra 3001-501, Portugal}
\address[3]{ Inria, Univ. Bordeaux, CNRS, Bordeaux INP, IMB, UMR 5251, 33405 Talence cedex, France}
\cortext[cor1]{Corresponding author (\href{mailto:mirco.ciallella@ensam.eu}{mirco.ciallella@ensam.eu})}

\begin{abstract}
In this paper we present a novel approach for the design of high order general boundary conditions 
when approximating solutions of the Euler equations on domains with curved boundaries, using meshes which may not be boundary conformal.
When dealing with curved boundaries and/or unfitted discretizations, the consistency of boundary conditions is a well-known challenge, especially 
in the context of high order schemes. In order to tackle such consistency problems, the so-called Reconstruction for Off-site Data (ROD) method has been recently 
introduced in the finite volume framework: it is based on performing a boundary polynomial reconstruction that embeds 
the considered boundary treatment thanks to the implementation of a constrained minimization problem. 
This work is devoted to the development of the ROD approach in the context of discontinuous finite elements.  
We use the genuine space-time nature  of the local ADER predictors to reformulate the ROD as a single space-time reconstruction procedure. 
This allows us to avoid a new reconstruction (linear system inversion) at each sub-time node and retrieve a single space-time polynomial
that embeds the considered boundary conditions for the entire space-time element.  
Several numerical experiments are presented proving the consistency of the new approach for all kinds of boundary conditions.
Computations involving the interaction of shocks  with  embedded curved boundaries are  made possible 
through an {\it a posteriori} limiting technique.
\end{abstract}

%% \begin{highlights}
%% \item Research highlight 1
%% \item Research highlight 2
%% \end{highlights}


\end{frontmatter}


\section{Introduction}
\label{sec:introduction}

The consistency of boundary conditions poses a significant challenge when dealing with high order discretization techniques. 
High-order methods have demonstrated their potential for providing more precise outcomes with fewer degrees of 
freedom~\cite{Wang_et_al_ijnmf13}. However, the management of boundaries remains an open challenging issue. 
In fact, when curved boundaries are taken into consideration, the precise representation of geometries and the accuracy of simulations become 
critical concerns. In particular, approximations of boundary conditions must be considered to achieve the desired convergence properties.

Currently, the most traditional method to address these issues employs high order polynomial 
reconstructions~\cite{zienkiewic} of the same degree used in 
the mathematical model discretization to estimate the targeted geometry. Additionally, other common approaches include local approximation 
of the curved geometry or the so-called iso-geometric analysis~\cite{nurbs0}. 
Curvilinear mesh generation presents a variety of challenges~\cite{luo2004automatic,sahni2010curved,doi:10.2514/6.2022-0389,doi:10.2514/6.2021-2991}, including nonlinear mappings associated with these elements and the 
use of complex quadrature formulas.  Methodological complexity, along with these features, 
substantially contributes to increase the computational costs in particular because of the augmented number of required degrees of freedom.

Despite the advancements made in developing numerical methods to handle curvilinear 
elements~\cite{dey2001towards,FORTUNATO20161,MOXEY2016130,Puigt}, obtaining a high-quality curved mesh 
remains a multifaceted task, particularly in dealing with practical geometric shapes. 
Alternatively, improving boundary conditions on a straight-faced mesh can be much cheaper, but it would need special care 
by taking into account the specific characteristics of the geometry in that area. 
We refer to~\cite{WangSun,krivodonova2006high} for the first works on curvature corrected wall boundary conditions 
and the application to high order schemes. 
In the previous reference, the boundary normal used when prescribing the slip-wall condition is correctly taken into account, 
resulting in a recovery of third and, in some cases, fourth order of accuracy on 2D geometries.                      
However, it is well-known that this method can only be formulated for slip-wall conditions, and the research is restricted to 2D geometries.  

More recent works on this research topic can be found in~\cite{ciallella2023shifted,costa2018very}.
In~\cite{ciallella2023shifted} some of the authors introduced a simple approach to design boundary conditions by reformulating the
Shifted Boundary Method (SBM)~\cite{Scovazzi1,Scovazzi2,Scovazzi3,lishifted,atallah2022high,carlier2023enriched,Scovazzi4,ciallella2022extrapolated,assonitis2022extrapolated} as a polynomial correction.
Although its simplicity to deal with Dirichlet conditions, the development of general Robin boundary conditions has only been 
addressed recently~\cite{visbech2023spectral}.
With the goal of considering general boundary conditions, some of the authors developed the Reconstruction for Off-site Data (ROD) 
approach in a high order finite volume framework~\cite{costa2018very}.
The ROD approach replaces the polynomial reconstruction in each boundary cell with a constrained least squares problem. 
The minimization problem  extracts information from neighboring cells and imposes constraints tied to the physical boundary. 
The work has been extended to other types of boundary conditions due to its generality~\cite{costa2019very,fernandez2020very,clain2021very,costa2021efficient,costa2022very,costa2023imposing}. 
In the finite volume framework, the constrained optimization problem is constructed by employing, as constraints, the information provided by the cell averages of the neighboring elements.  

In this paper, we investigate the development of the ROD approach within a discontinuous finite element framework~\cite{cockburn1991runge,cockburn1998runge,cockburn2001runge,cockburn2012discontinuous,mazaheri2019bounded}. 
Specifically, by having the local formulation of discontinuous finite element discretizations, a distinct constrained minimization 
problem can be formulated which only employs the information of the element under consideration. 
This simplifies certain aspects of the finite volume approach while maintaining its consistency. 
To achieve high order properties for time-dependent problems, the approach based on spatial polynomial 
reconstruction requires the solution of this minimization problem at every sub-time step (when coupled with Runge-Kutta schemes) 
of the considered time integration scheme, resulting in increased computational costs.
To further simplify this process, a novel approach is introduced in the space-time ADER framework~\cite{dumbser2008unified,boscheri2014high,boscheri2017high,boscheri2015direct,busto2020high,fambri2017space,balsara2009efficient,gaburro2020high}. 
In particular, we propose a one-step direct and cheap procedure by employing the space-time basis functions used in the ADER predictor 
also to write the ROD constrained least squares minimization problem.
The obtained system requires only one inversion and provides a single polynomial for the entire space-time element.
Due to the modal  space-time representation used, the number of modes reconstructed is much  reduced compared to  the case of several independent reconstructions at each quadrature point in time. 
This provides a cost effective formulation of ROD  for evolutionary problems. 
Numerical results confirm the accuracy of the proposed approach on unfitted meshes.
 
The paper is organized as follows. In~\cref{sec:model}, we present the studied mathematical model, 
specifically the hyperbolic system of the Euler equations of gasdynamics. 
However, it should be noted that the proposed boundary treatment can be suitably applied to other PDEs. 
In~\cref{sec:numerics}, we review the high order discretization of conservation laws system using 
discontinuous finite element~\cite{cockburn1998runge} with modal polynomials and the one-step ADER method 
with space-time predictor~\cite{dumbser2008unified}.  
A shock limiting approach~\cite{clain2011high,diot2013multidimensional,loubere2014new,gaburro2021posteriori} 
is also discussed in this section to deal with discontinuous features that may arise in the flow field.
Additionally, in~\cref{sec:boundarynumerics}, we introduce the boundary treatment based on the 
Reconstruction for Off-site Data (ROD) technique. We begin by examining the classical approach adapted to 
the discontinuous Galerkin discretization, which has not yet been explored in the literature. 
Then, we introduce a new method that utilizes the space-time high order predictor to create a boundary treatment 
that eliminates the need for polynomial reconstruction at each sub-time node of the ADER correction step. 
In~\cref{sec:results}, we present numerical experiments demonstrating the convergence properties of the methods on 
conformal linear meshes and unfitted triangulations. 
The paper ends with~\cref{sec:conclusions}, which offers an outlook regarding the future 
analysis and development of these schemes.

\section{Mathematical model}
\label{sec:model}

The mathematical model considered in this work is the system of Euler equations for compressible gasdynamics in $d=2$ space dimensions reading:
\begin{equation}\label{eq:euler0}
\frac{\partial\, \mathcal{U}}{\partial t}+\nabla\cdot\mathcal{F}(\mathcal{U}) =0, \quad \text{on}\quad \Omega_T=\Omega\times[0,T]\subset\mathbb{R}^d\times\mathbb{R}^+,
\end{equation}
with  $\mathcal{U}:\R^+\times\R^d\rightarrow\R^{d+2}$ the vector of conserved variables and $\mathcal{F}:\R^{d+2}\rightarrow\R^d\times\R^{d+2}$ the nonlinear flux, respectively defined as
\begin{equation}\label{eq:euler0a}
\mathcal{U}=\left(
\begin{array}{c}
\rho \\ \rho \mathbf{u}  \\ \rho E
\end{array}
\right)\;,\;\;
\mathcal{F}(\mathcal{U})= \left(
\begin{array}{c}
\rho \mathbf{u}\\ \rho \mathbf{u} \otimes \mathbf{u} + p \mathbb{I} \\ \rho H \mathbf{u} 
\end{array}
\right), 
\end{equation}
having denoted by $\rho$ the mass density, by $\mathbf{u}$ the velocity, by $p$ the pressure, and with $E=e+\mathbf{u}\cdot \mathbf{u}/2$ the specific total energy, $e$ being the specific internal energy.
Finally, the total specific enthalpy is  $H=h+\mathbf{u}\cdot \mathbf{u}/2$ with  $h=e+p/\rho$ the specific enthalpy. 
The relation between the pressure and the internal energy is given by the perfect gas equation of state%:
\begin{equation}\label{eq:EOS}
p=(\gamma-1)\rho e,
\end{equation}
where $\gamma$ the constant ratio of specific heats ($\gamma=1.4$ for air).\\

The validation of the numerical methods are performed by means of the manufactured solution method which needs the discretization of an additional source term 
in~\cref{eq:euler0} which thus reads as follows
\begin{equation}\label{eq:euler1}
 \frac{\partial\,\mathcal{U}}{\partial t}+\nabla\cdot\mathcal{F}\left(\mathcal{U}\right) = \mathcal{S}(\mathcal{U}).
\end{equation}
For this reason, we are going to include $\mathcal{S}$ in the space-time discretization presented in~\cref{sec:numerics}.

\section{High order space-time discretization}
\label{sec:numerics}

\subsection{Computational domain and data representation}

We discretize the domain $\Omega$ as a tessellation $\mathscr{T}$ composed of $\mathscr{N}$ non-overlapping simplicial  elements (triangles in 2D).
We denote by $K_j$ the generic element, so that the computational domain is $\Omega_h=\bigcup_{j=1}^{\mathscr{N}} K_j$. 
Note that  in general  $\Omega_h\ne \Omega$ and in particular $\partial\Omega_h\ne\partial\Omega$ for most approximations, even conformal, 
except for very simple geometries (with no curvature) or of iso-geometric approaches~\cite{nurbs0}.

The solution $\mathcal{U}$ is approximated by $\mathcal{U}_h$, which belongs to a space of piece-wise polynomials within each triangle $K_j$ and discontinuous across faces, such that in each element we have
\begin{equation}\label{eq:basis psi}
\mathcal{U}_h(\mathbf{x},t)|_{K_j} = \sum_{k}\psi_{k}(\mathbf{x}) \,\hat{\mathcal{U}}_{k}(t), 
\end{equation}
where $\psi_{k}$ is a basis of polynomials of degree $M$.
It is well-known that the discontinuous finite element data representation~\eqref{eq:basis psi} leads to a Finite Volume (FV) scheme if $M=0$.
\subsection{Discontinuous Galerkin in space}
The elemental semi-discrete discontinuous Galerkin (DG) weak formulation is classically written by projecting each component of~\cref{eq:euler1}
on the relevant basis and integrating by parts~\cite{cockburn1998runge,bassi1997high}
\begin{equation}\label{eq:DGweak semidiscrete 2}
\int_{K_j} \psi \, \frac{\diff{\,\mathcal{U}_h}}{\diff{t}} \diff{\mathbf{x}} + \int_{\partial K_j} \psi \, \hat{F}(\mathcal{U}_h^-,\mathcal{U}_h^+)\cdot \mathbf{n}\diff{S}-\int_{K_j} \nabla\psi \, \cdot\mathcal{F}(\mathcal{U}_h)\diff{\mathbf{x}} = \int_{K_j} \psi \, \mathcal{S}(\mathcal{U}_h)\diff{\mathbf{x}}, 
\end{equation}
with $\hat{\mathcal F}(\mathcal{U}_h^-,\mathcal{U}_h^+)$ a consistent numerical flux which depends on the face values of the internal  state $\mathcal{U}_h^-$, 
on the neighboring element  state $\mathcal{U}_h^+$ and on the face normal $\mathbf{n}$. 
We recall in particular that a consistent flux  is a Lipschitz continuous function whose arguments also verify
\begin{equation}\label{eq:constistencyDG}
\hat{\mathcal F}(U,U)\cdot \mathbf{n} = \mathcal{F}(U)\cdot \mathbf{n}.
\end{equation}
In this paper, we have used a classical Rusanov-type (local Lax-Friedrichs) flux:
\begin{equation}\label{eq:rusanov flux}
\hat{\mathcal F}(\mathcal{U}_h^-,\mathcal{U}_h^+)\cdot \mathbf{n} = \frac{1}{2}\left(\mathcal{F}(\mathcal{U}_h^+) + \mathcal{F}(\mathcal{U}_h^-)\right)\cdot \mathbf{n} - \frac{1}{2}s_{max}\left(\mathcal{U}_h^+ -\mathcal{U}_h^-\right),
\end{equation}
where $s_{max}$ is the largest of the spectral radiuses of the flux Jacobians $\frac{\partial\mathcal F}{\partial\mathcal U} (\mathcal{U}_h^+)$ and
$\frac{\partial\mathcal F}{\partial\mathcal U} (\mathcal{U}_h^-)$. %\\
\subsection{ADER in time with local space-time predictor}
For the time integration, we have used a high order explicit predictor-corrector ADER method~\cite{dumbser2008unified}.
The time domain is discretized in temporal intervals  $[t^n, t^{n+1}]$, where $t^{n+1}=t^n+\Delta t$.
%The time step size $\Delta t$ is chosen with the classical explicit CFL condition~\eqref{eqn.cfl}. 
With the same notation used before, the ADER method can be recast as
\begin{equation}\label{eq:corrector2}
  \mathcal{U}^{n+1}  = \mathcal{U}^{n} -  M^{-1} \int_{t^n}^{t^{n+1}} R(q_h) dt, % 
\end{equation}
where $M$ is the elemental mass matrix, $R$ includes all integrals in~\cref{eq:DGweak semidiscrete 2} but the first, and $q_h=q_h(\mathbf{x},t)$ is  a high order space-time predictor of the solution valid in the time interval $[t^n, t^{n+1}]$.
In practice, \cref{eq:corrector2} is replaced by a consistent high order quadrature formula in time (with $\alpha_i$ and $\omega_i$, the quadrature points and weights)
\begin{equation}\label{eq:corrector2a}
  \mathcal{U}^{n+1}  =    \mathcal{U}^{n} -  M^{-1}  \Delta t 
  \sum\limits_{i}\omega_i R\big(q_h(t^{n+\alpha_i}) \big). % 
\end{equation}
$q_h(\mathbf{x},t)$ is defined as a polynomial of degree $M$ in space and time (where time is considered as an additional physical dimension),
and $R$ also includes consistent quadrature formulas for the integrals in space.
This polynomial is obtained by means of a genuinely local space-time procedure.
In the current implementation this local problem  is formulated  by means of a modal expansion in the space-time element
\begin{equation}\label{eq:predictor0}
q_h(\mathbf{x},t)|_{K_j\times[t^n, t^{n+1}]} = \sum_{\ell=0}^{\mathcal{Q}-1}\theta_\ell(\mathbf{x},t)\hat{q}_\ell, \qquad \mathcal{Q}=\mathcal{L}(M,d+1),
\end{equation}
where  $\theta_\ell(\mathbf{x},t)$ are $\mathcal{L}(M,d)=\frac{1}{d!}\prod_{m=1}^d (M + m)$ modal bases defined as
\begin{equation}
\theta_\ell(\mathbf{x},t)_{K_j} = \frac{(x_1-x_{1,j}^n)^{p_\ell}}{p_\ell ! h_j^{p_\ell}} \frac{(x_2-x_{2,j}^n)^{q_\ell}}{q_\ell ! h_j^{q_\ell}} \frac{(t-t^n)^{r_\ell}}{r_\ell ! h_j^{r_\ell}}, \quad \ell=0,\ldots,\mathcal{L}(M,d+1),\quad 0\leq p_\ell+q_\ell+r_\ell\leq M,
\end{equation}
with $(x_{1,j},x_{2,j})$ the barycenter of element $K_j$, and $h_j$ is the diameter of its incircle.
 
The values $\hat q_\ell$ are obtained as the solution of the local space-time weak formulation
\begin{align}
\int_{K_j} \theta_k(\mathbf{x},t^{n+1}) \, q_h(\mathbf{x},t^{n+1})\diff{\mathbf{x}} &- \int_{K_j} \theta_k(\mathbf{x},t^{n})  \, \mathcal{U}^n_h(\mathbf{x},t^{n})\diff{x} - \int_{t^n}^{t^{n+1}}\int_{K_j} \frac{\partial \theta_k}{\partial t}(\mathbf{x},t) \, q_h(\mathbf{x},t)\diff{\mathbf{x}}\diff{t}  \nonumber\\ &+\int_{t^n}^{t^{n+1}}\int_{K_j} \theta_k(\mathbf{x},t) \, \nabla\cdot\mathcal{F}(q^n_h)\diff{\mathbf{x}}\diff{t}=\int_{t^n}^{t^{n+1}}\int_{K_j} \theta_k(\mathbf{x},t) \, \mathcal{S}(q^n_h)\diff{\mathbf{x}}\diff{t} , \label{eq:predictor2}
\end{align}
where $\mathcal{U}^n_h$ is the known initial condition at time $t^n$.

Equation~\eqref{eq:predictor2} is local, meaning that no communication is exchanged between $K_j$ and its neighboring control volumes. 
The solution to this equation can be obtained independently within each element $K_j$ by means of some iterative procedure.
Herein, a simple discrete Picard iteration for each space-time element is performed~\cite{busto2020high}.

\subsection{\textit{A posteriori} sub-cell finite volume limiter}\label{sec:limiter}

In order to properly handle discontinuous solutions we employ the MOOD approach~\cite{clain2011high,diot2012improved,diot2013multidimensional},
which has already been effectively applied in the ADER framework~\cite{loubere2014new,boscheri2015direct,boscheri2017high,gaburro2021posteriori}.

The algorithm is based on an \textit{a posteriori} technique.
The solution is first evolved from $t^n$ to $t^{n+1}$ using the high order ADER-DG method. 
Then several admissibility criteria are checked and the solution in all troubled cells (i.e., the cells not satisfying such criteria) 
is recomputed  \textit{a posteriori} using a robust second-order MUSCL-Hancock TVD finite volume scheme, 
on a sub-triangulation of the initial grid.
This allows us to preserve the accuracy of the high order DG scheme also when passing to a lower order but more robust scheme.
All aspects of the implementation of this technique are provided in~\cite{gaburro2021posteriori} to which we refer for details. \\
\begin{remark}[CFL constraint]
Concerning the choice of the time step, we have implemented the classical explicit CFL condition 
\begin{equation}
\label{eqn.cfl}
\Delta t < \text{CFL} \left( \frac{|h_{\text{min}}|}{(2M+1)|\lambda_{\text{max}}|}  \right),
\end{equation}
where $|h_{min}|$ is the minimum characteristic mesh-size and $|\lambda_{max}|$ is the maximum eigenvalue of the Jacobian of the flux.
For DG on unstructured meshes the CFL stability condition requires a $CFL<1/d$ to be satisfied (we refer to~\cite{CHALMERS2020109095} for more details).
Notice that the time step constraint does not need to be modified in presence of troubled cells, 
because we subdivide each troubled triangle in exactly $(2M+1)^d$ sub-triangles, and then we employ a FV scheme for which~\eqref{eqn.cfl} holds with $M=0$.
\end{remark}

\section{Boundary conditions}

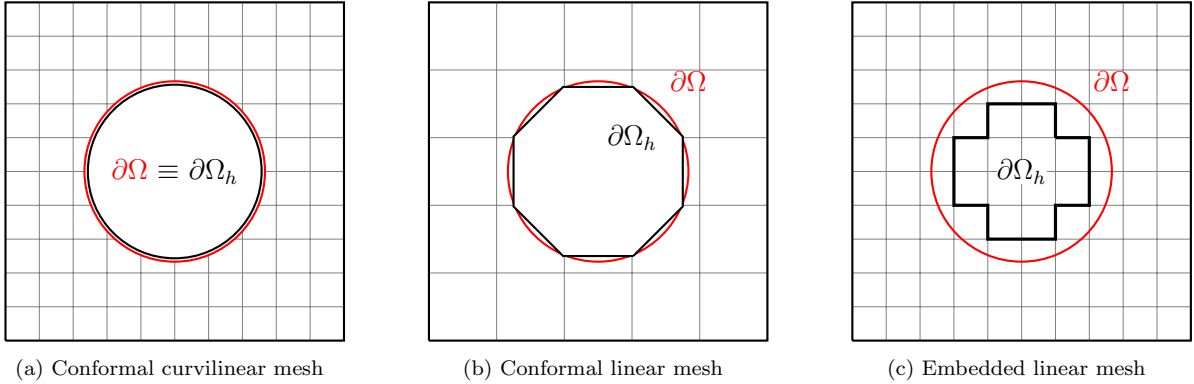
\begin{figure}
\centering
\subfloat[Conformal curvilinear mesh]{
        \begin{tikzpicture}[scale=1.5]
 		\draw[gray,very thin,step=0.3] (0,0) grid (3,3);
		\draw[thick] (0.,0) -- (0,3) -- (3,3) -- (3,0) -- (0,0);
		\draw[thick,red,fill=white] (1.5,1.5) circle (0.8cm);
		\draw[thick] (1.5,1.5) circle (0.77cm);
	        \node[scale=0.9] at (1.5,1.5) {\color{red}{$\partial\Omega$} \color{black}{$\equiv$} \color{black}{$\partial\Omega_h$}};
        \end{tikzpicture}
}\qquad
\subfloat[Conformal linear mesh]{
        \begin{tikzpicture}[scale=1.5]
	        %\node [black] at (0.5,0.5) {$\Omega$};
 		\draw[gray,very thin,step=0.6] (0,0) grid (3,3);
		\draw[thick] (0.,0) -- (0,3) -- (3,3) -- (3,0) -- (0,0);
		\draw[thick,red] (1.5,1.5) circle (0.8cm);
		\draw[thick,fill=white] (1.5,1.5) \polygon{8}{0.75};
	        \node [black,scale=0.9] at (1.8,1.8) {$\partial\Omega_h$};
	        \node [red,fill=white,scale=0.9] at (2.3,2.3) {$\partial\Omega$};
        \end{tikzpicture}
}\qquad
\subfloat[Embedded linear mesh]{
	\begin{tikzpicture}[scale=1.5]
 		\draw[gray,very thin,step=0.3] (0,0) grid (3,3);
		\draw[thick] (0.,0) -- (0,3) -- (3,3) -- (3,0) -- (0,0);
		\draw[thick,red] (1.5,1.5) circle (0.8cm);
	        \node [red,fill=white,scale=0.9] at (2.3,2.3) {$\partial\Omega$};
		\draw[very thick] (1.8,1.8) -- (2.1,1.8) -- (2.1,1.2) -- (1.8,1.2) -- (1.8,0.9) -- (1.2,0.9) -- (1.2,1.2) -- (0.9,1.2) -- (0.9,1.8) -- (1.2,1.8) -- (1.2,2.1) -- (1.8,2.1) -- (1.8,1.8); 
	        \node [black,fill=white,scale=0.9] at (1.5,1.5) {$\partial\Omega_h$};
	\end{tikzpicture}
}
\caption{The different possibilities available for the discretizations of domains with internal boundaries.}\label{fig:meshbndexample}
\end{figure}

In this section, we firstly discuss the classical approach to enforce boundary conditions in a DG framework.
In particular, we focus on Dirichlet-type and slip-wall conditions.
However, this technique is consistent as long as the boundary is also discretized with high order polynomials (see~\cref{fig:meshbndexample}a).
When discretizing it with a linear conformal mesh, as in~\cref{fig:meshbndexample}b, this approach is at best second order accurate.
Finally, by using fully embedded meshes, as in~\cref{fig:meshbndexample}c, the boundary consistency automatically drops to first order if no special treatments are considered.   
We then consider a more general slip-wall condition to design high order consistent boundary treatments for the cases described in~\cref{fig:meshbndexample}b--c. 

\subsection{Classical boundary treatment}

When the boundary of element $K_j$ belongs to $\partial\Omega_h$, the normal flux function must account for the appropriate  boundary conditions. 
In this work, the boundary flux is obtained by defining a ghost state $\mathcal{U}^{bc}$,
 and introducing a numerical flux $\hat{\mathcal F}(\mathcal{U}^-_h,\mathcal{U}^{bc})$ defined by an approximate Riemann solver~\eqref{eq:rusanov flux} based 
on the internal state $\mathcal{U}^-_h$ and on $\mathcal{U}^{bc}$. 
Depending on the condition to be enforced, different definitions of the ghost state are used.
In this work we focus on  
\begin{itemize}
\item Dirichlet-type conditions, where $\mathcal{U}^{bc}$ is fully known and enforced weakly through fluxes;
\item reflecting conditions (slip-wall), where we set $\mathbf{u}\cdot \mathbf{n} = 0$, for static walls.
In general for conformal meshes, $\mathcal{U}^{bc}$ has the same density, internal energy and tangential velocity of $\mathcal{U}^-_h$, and 
the opposite normal velocity component $\mathbf{u}\cdot \mathbf{n}$. 
\end{itemize}
For conformal meshes, it is then consistent to take simply the face normal $\mathbf{n}$ as the normal of $\partial\Omega_h$ pointing outwards.

%%  \begin{remark}[Boundary numerical flux]
%%  In particular, we find beneficial for very high-order schemes to use a variation of the Rusanov flux that does not include the ghost state within
%%  the flux average, meaning
%%  %
%%  \begin{equation}
%%  \hat{\mathcal F}(\mathcal{U}^-_h,\mathcal{U}^{bc}) = \mathcal F(\mathcal{U}^-_h)\cdot n - \frac12 s_{max} (\mathcal U^{bc}_h - \mathcal U^-_h). 
%%  \end{equation}
%%  %
%%  We believe that by removing the reconstructed ghost state from the nonlinearity of the flux average, the conditioning of the system related
%%  to the numerical boundary conditions greatly improves. 
%%  \end{remark}
%%  %
%%  \begin{verbatim}
%%  Maybe remove the last remark
%%  \end{verbatim}

\subsection{General slip-wall boundary treatment}

When working on conformal curvilinear meshes, the reflecting boundary condition can be consistently formulated by imposing $\mathbf{u}\cdot\mathbf{n} = 0$
thanks to the high order approximation of the boundary (see~\cref{fig:meshbndexample}a).
However, when approximating linearly or with general embedded meshes curved boundaries the condition imposed on the normal component 
of the velocity might not be enough due to the distance between $\partial\Omega$ and $\partial\Omega_h$ 
(see~\cref{fig:meshbndexample}b--c). 
For this reason, in this work we consider the more general boundary conditions introduced in~\cite{chertock2018second} 
to have an accurate slip-wall treatment.
Herein we recall the basic steps and hypothesis.

Starting from the mass and momentum in~\cref{eq:euler0}, under some smoothness assumptions on the conserved variables we can write,

\begin{equation}
\rho\frac{D\mathbf{u}}{Dt} + \nabla p = 0,
\end{equation}
where $\frac{D}{Dt} = \frac{\partial}{\partial t}  + (\mathbf{u}\cdot\nabla)$ represents the material derivative.

The first condition to be imposed is on the normal component of the velocity vector for reflecting boundary conditions, meaning 

\begin{equation}\label{eq:slip1}
u_n = \mathbf{u}\cdot\mathbf{n} = 0,
\end{equation}
where $\mathbf{n}$ is now the normal of the real geometry, which may greatly differ from the normal of the discretized boundary $\partial\Omega_h$.

It should be noticed that in this work the real boundary normal is defined starting from a level set $\phi$ initialized as a signed distance function,
$$ \mathbf{n}:=\frac{\nabla\phi}{\|\nabla\phi\|},$$
and the real tangent vector is defined such that $\mathbf{n}\cdot \boldsymbol{\tau} = 0$, and $\boldsymbol{\tau}\times \mathbf{n} = 1$.

Condition~\eqref{eq:slip1} implies that we want to enforce the flow to be tangential to the real boundary $\partial\Omega$,
meaning $\mathbf{u} = u_\tau \boldsymbol{\tau}$ which allows us to write,

\begin{equation}
\frac{D\mathbf{u}}{Dt} = \boldsymbol{\tau} \frac{\partial u_\tau}{\partial t} + u_\tau \frac{\partial \boldsymbol{\tau}}{\partial t} + u_\tau (\mathbf{u}\cdot\nabla) \boldsymbol{\tau} = a_\tau \boldsymbol{\tau} + u_\tau \frac{D\boldsymbol{\tau}}{Dt}, 
\end{equation}
where $a_\tau$ is the tangential acceleration.
By using the Frenet-Serret formulae $\frac{D\boldsymbol{\tau}}{Dt} = u_\tau k \mathbf{n}$ (where $k$ is the local curvature), 
and projecting the momentum equation onto the normal direction, we end up with the second condition on the pressure,

\begin{equation}\label{eq:slip2}
\frac{\partial p}{\partial n} = -\rho u_\tau^2 k.
\end{equation}
%
%% \begin{remark}[Material derivative of the tangent]
%% The relation on the material derivative of the tangent vector is recovered considering the Frenet-Serret formulae written for a parametrized curve $r(s)$.
%% Starting from the definition of a tangent unit vector $r'(s)=\tau(s)$, we can write $\|\tau(s)\|=\tau(s)\cdot\tau(s)=1$.
%% By taking the derivative of the scalar product we obtain $\tau'(s)\cdot \tau(s) = 0$, which entails that $\tau'(s)$ is parallel to the normal $n(s)$.
%% From the normalized relation we deduce the Frenet-Serret formula $\tau'(s) = \|\tau'(s)\| n(s) = k(s) n(s) $.
%% Finally, by differentiating with respect to time the tangent vector we obtain
%% $$ \frac{\mathrm{d}}{\mathrm{d}t}\tau(s(t)) = \frac{\mathrm{d}\tau}{\mathrm{d}s} \frac{\mathrm{d}s}{\mathrm{d}t} = \tau'(s) u(s) =  k(s) n(s) u(s)  $$ 
%% \end{remark}
%
The third condition on the density is found by imposing an adiabatic wall, meaning $\partial_n S=0$ where S denotes the entropy.
From $S=p\rho^{-\gamma}$ we can easily obtain $\partial_n S = (\partial_n p) \rho^{-\gamma} - \gamma p \rho^{-\gamma-1} (\partial_n \rho)$, and

\begin{equation}\label{eq:slip3}
\frac{\partial \rho}{\partial n} = \frac1{c^2} \frac{\partial p}{\partial n},
\end{equation}
 
where $c = \sqrt{\gamma p/\rho}$ is the speed of sound.

In order to derive the condition on the tangential component of the velocity, 
we consider the vorticity $\boldsymbol\omega:=\nabla\times \mathbf{u}$ 
and compute it in the local coordinates $(\mathbf{n},\,\boldsymbol{\tau})$ considering $\mathbf{u} = u_n \mathbf{n} + u_\tau \boldsymbol{\tau}$, 

%% $$ \boldsymbol\omega = \frac{\partial u_\tau}{\partial n} - \frac{\partial u_n}{\partial \tau}  + u_n \nabla\times n + u_\tau \mathrm{div}(n) = 
%% \frac{\partial u_\tau}{\partial n} - \frac{\partial u_n}{\partial \tau} -u_\tau k, $$
%% where, in the last step, we used the definition of the normal $\nabla\times(\nabla\phi)=0$, and $k=-\mathrm{div}(n)$.

$$ \boldsymbol\omega = \frac{\partial u_\tau}{\partial n} - \frac{\partial u_n}{\partial \tau} -u_\tau k. $$

Finally, by assuming that the flow is irrotational and $u_n=0$, we can deduce 
the needed condition on $u_\tau$,

\begin{equation}\label{eq:slip4}
\frac{\partial u_\tau}{\partial n} = u_\tau k.
\end{equation}
  
In~\cite{chertock2018second} it was also shown that the condition on $u_\tau$ can be retrieved by imposing $\partial_n H = 0$.
\begin{remark}[Equivalent boundary conditions]
From the previous statements, it is straightforward to deduce that the same boundary conditions can be imposed through the more 
compact system of mixed conditions:
\begin{equation}\label{eq:slipeqv}
u_n = 0,\quad \frac{\partial S}{\partial n} = 0, \quad \frac{\partial H}{\partial n} = 0,\quad \frac{\partial u_\tau}{\partial n} - u_\tau k = 0.
\end{equation}
Notice that by performing a change of variables from the conservative ones  to $(u_n,\,S,\,H,\,u_\tau)$ we can use the mixed Dirichlet-Neumann-Robin conditions~\eqref{eq:slipeqv} 
to compute the ghost state and impose the slip-wall boundary condition.
\end{remark}
\section{High order ROD boundary reconstruction}
\label{sec:boundarynumerics}

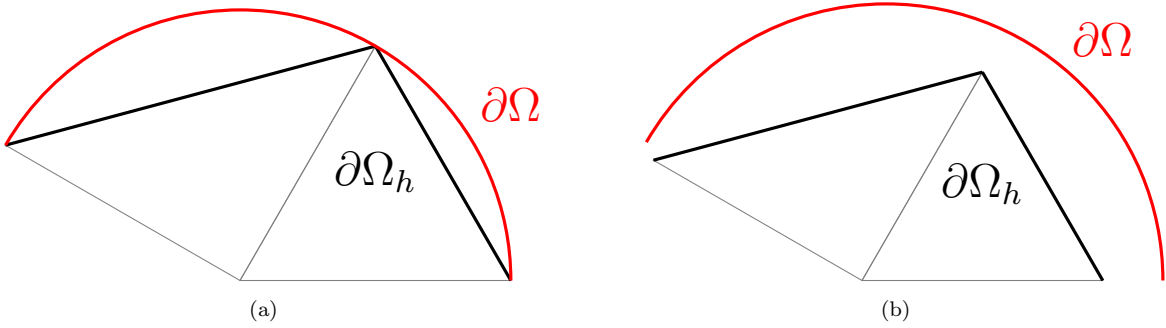
\begin{figure}
\centering
\subfloat[]{
\begin{tikzpicture}[scale=1.8]
  % Definizione dei vertici dei triangoli
  \coordinate (A) at (0,0);
  \coordinate (B) at (2,0);
  \coordinate (C) at (1,{sqrt(3)});
  \coordinate (D) at ({-sqrt(3)},1);

  % Disegno dei triangoli
  \draw[gray] (A) -- (B) -- (C) -- cycle;
  \draw[gray] (A) -- (C) -- (D) -- cycle;
  \draw[very thick] (B) -- (C);
  \draw[very thick] (D) -- (C);

  % Curva che tocca i vertici dei triangoli
  \draw[red, very thick] (B) arc (0:150:2); 
  \node[red,fill=white,scale=1.5] at (2.0,1.3) {$\partial\Omega$};
  \node[fill=white,scale=1.5] at (1,0.8) {$\partial\Omega_h$};

\end{tikzpicture}}\qquad\quad
\subfloat[]{
\begin{tikzpicture}[scale=1.6]
  % Definizione dei vertici dei triangoli
  \coordinate (A) at (0,0);
  \coordinate (B) at (2,0);
  \coordinate (C) at (1,{sqrt(3)});
  \coordinate (D) at ({-sqrt(3)},1);

  % Disegno dei triangoli
  \draw[gray] (A) -- (B) -- (C) -- cycle;
  \draw[gray] (A) -- (C) -- (D) -- cycle;
  \draw[very thick] (B) -- (C);
  \draw[very thick] (D) -- (C);

  % Curva che tocca i vertici dei triangoli
  \draw[red, very thick] (2.5,0) arc (0:150:2.3); 
  \node[red,fill=white,scale=1.5] at (2.0,2) {$\partial\Omega$};
  \node[fill=white,scale=1.5] at (1,0.8) {$\partial\Omega_h$};

\end{tikzpicture}}
\caption{Visual representation of a curved boundary discretized with a linear conformal triangulation (a) and a fully embedded (b) triangulation.}\label{fig:confembtriangexample}
\end{figure}

In this section, we present the numerical treatments we propose herein to impose general boundary conditions.
We begin by presenting a reformulation of the ROD boundary reconstruction in space by employing a set of traditional basis function, which differ from those used for the discretization of the spatial operators.  
Afterwards, we introduce a novel space-time ROD reconstruction that exploits the space-time basis function used to define the high order local predictor~\eqref{eq:predictor0}.
It is then shown that performing the ROD reconstruction by using the same basis functions used for the discretization allows us to greatly simplify the functional of the ROD minimization problem.

%It should be noticed, that no differences in the ROD algorithm are needed when dealing with linear conformal and fully embedded 
%triangulations (see~\cref{fig:confembtriangexample}): the only difference is the larger distance of the considered quadrature point
%with respect to the real geometry, which is something that both approaches always take into account.

\subsection{Reconstruction for off-site data methods in space}

The Reconstruction for Off-site Data (ROD) method consists in discarding the polynomial of the boundary cell
$K_j\in \partial\Omega_h$ and replacing it with a new reconstructed one that embeds the considered boundary condition.
The ROD polynomial is then used to reconstruct the solution in each quadrature point $\mathbf{x}^\star$ of the boundary integral on $\partial\Omega_h$. 
This reconstruction is performed through a constrained least-squares procedure, which is then solved with Lagrange multipliers.  

\begin{figure}
\centering
\subfloat[ROD at each sub-time node ($\textcolor{blue}{\rightarrow}$)]{
	\begin{tikzpicture}[scale=2.5]
	\draw[thick,-stealth,-dot-=0,-dot-=0.83,-dot-=0.3,-dot-=0.6] (-1,0) -- (-1,2.4);
	\draw[thick,->,blue] (-1.8,0)  -- (-1.5,0) ;
	\draw[thick,->,blue] (-1.8,0.7)  -- (-1.5,0.7) ;
	\draw[thick,->,blue] (-1.8,1.4)  -- (-1.5,1.4) ;
	\draw[thick,->,blue] (-1.8,2.0)  -- (-1.5,2.0) ;
	\node[] at (-1.25,0) {$t^n$};
	\node[] at (-1.25,0.7) {$t^{n+\alpha_1}$};
	\node[] at (-1.25,1.4) {$t^{n+\alpha_2}$};
	\node[] at (-1.25,2) {$t^{n+1}$};
	%\node[blue] at (-1.3,1.1) {$t_{gp}$};
	%\node[red] at (0.7,2.4) {$\Gamma$};
	%\draw[blue,decorate, decoration = {brace}] (-1.05,0.71) --  (-1.05,1.42);
	% arc
	\draw[thick,red,-dot-=0.5] (-0.75,2) to[bend left] (0,2.5);
	\draw[thick,red,-dot-=0.5] (0,2.5) to[bend left] (0.75,2);
	\draw[dashed,thick,red,-dot-=0.5] (-0.75,0) to[bend left] (0,0.5);
	\draw[dashed,thick,red,-dot-=0.5] (0,0.5) to[bend left] (0.75,0);
	\draw[dashed,thick,red,-dot-=0.5] (-0.75,0.7) to[bend left] (0,1.2);
	\draw[dashed,thick,red,-dot-=0.5] (0,1.2) to[bend left] (0.75,0.7);
	\draw[dashed,thick,red,-dot-=0.5] (-0.75,1.4) to[bend left] (0,1.9);
	\draw[dashed,thick,red,-dot-=0.5] (0,1.9) to[bend left] (0.75,1.4);
	% stages
	\draw[thick,dashdotted,blue] (-0.75,0.7) -- (0,1.2) -- (0.75,0.7)  -- cycle;
	\draw[thick,dashdotted,blue] (-0.75,1.4) -- (0,1.9) -- (0.75,1.4)  -- cycle;
	
	\foreach \Point in {(-0.4,2.1), (0,2.1), (0.4,2.1), (-0.2,2.25) , (0,2.25),  (0.2,2.25) , (0,2.4) }{
	    \node[blue] at \Point {\textbullet};
	}

	\foreach \Point in {(-0.4,1.5), (0,1.5), (0.4,1.5), (-0.2,1.65) , (0,1.65),  (0.2,1.65) , (0,1.8) }{
	    \node[blue] at \Point {\textbullet};
	}
	
	\foreach \Point in {(-0.4,0.8), (0,0.8), (0.4,0.8), (-0.2,0.95) , (0,0.95),  (0.2,0.95) , (0,1.1) }{
	    \node[blue] at \Point {\textbullet};
	}
	
	\foreach \Point in {(-0.4,0.1), (0,0.1), (0.4,0.1), (-0.2,0.25) , (0,0.25),  (0.2,0.25) , (0,0.4) }{
	    \node[blue] at \Point {\textbullet};
	}
	% triangular prism
	\draw[dashed] (-0.75,0) -- (0,0.5) edge (0,2.5) -- (0.75,0);
	\draw[thick] (-0.75,0) rectangle (0.75,2) -- (0,2.5) -- (-0.75,2);
	\end{tikzpicture}}\qquad\qquad
\subfloat[ROD in the space-time element ($\textcolor{darkgreen}{\rightarrow}$)]{
	\begin{tikzpicture}[scale=2.5]
	\draw[thick,-stealth,-dot-=0,-dot-=0.83,-dot-=0.3,-dot-=0.6] (-1,0) -- (-1,2.4);
	\node[] at (-1.25,0) {$t^n$};
	\node[] at (-1.25,2) {$t^{n+1}$};
	\draw[thick,->,darkgreen] (-1.9,1.0)  -- (-1.6,1.0) ;
	%\node[red] at (0.7,2.4) {$\Gamma$};
	\draw[darkgreen,decorate, decoration = {brace}] (-1.5,0.0) --  (-1.5,2.0);
	% arc
	\draw[thick,red,-dot-=0.5] (-0.75,2) to[bend left] (0,2.5);
	\draw[thick,red,-dot-=0.5] (0,2.5) to[bend left] (0.75,2);
	\draw[dashed,thick,red,-dot-=0.5] (-0.75,0) to[bend left] (0,0.5);
	\draw[dashed,thick,red,-dot-=0.5] (0,0.5) to[bend left] (0.75,0);
	\draw[dashed,thick,red,-dot-=0.5] (-0.75,0.7) to[bend left] (0,1.2);
	\draw[dashed,thick,red,-dot-=0.5] (0,1.2) to[bend left] (0.75,0.7);
	\draw[dashed,thick,red,-dot-=0.5] (-0.75,1.4) to[bend left] (0,1.9);
	\draw[dashed,thick,red,-dot-=0.5] (0,1.9) to[bend left] (0.75,1.4);
	% stages
	\draw[thick,dashdotted,darkgreen] (-0.75,0.7) -- (0,1.2) -- (0.75,0.7)  -- cycle;
	\draw[thick,dashdotted,darkgreen] (-0.75,1.4) -- (0,1.9) -- (0.75,1.4)  -- cycle;
	
	%\foreach \Point in {(-0.4,2.1), (0,2.1), (0.4,2.1), (-0.2,2.25) , (0,2.25),  (0.2,2.25) , (0,2.4) }{
	\foreach \Point in {(-0.2,2.25) , (0,2.25),  (0.2,2.25) }{
	    \node[darkgreen] at \Point {\textbullet};
	}
	
	%\foreach \Point in {(-0.4,1.5), (0,1.5), (0.4,1.5), (-0.2,1.65) , (0,1.65),  (0.2,1.65) , (0,1.8) }{
	\foreach \Point in { (-0.2,1.65) , (0,1.65),  (0.2,1.65)  }{
	    \node[darkgreen] at \Point {\textbullet};
	}
	
	%\foreach \Point in {(-0.4,0.8), (0,0.8), (0.4,0.8), (-0.2,0.95) , (0,0.95),  (0.2,0.95) , (0,1.1) }{
	\foreach \Point in {(-0.2,0.95) , (0,0.95),  (0.2,0.95)  }{
	    \node[darkgreen] at \Point {\textbullet};
	}

        %\foreach \Point in {(-0.4,0.1), (0,0.1), (0.4,0.1), (-0.2,0.25) , (0,0.25),  (0.2,0.25) , (0,0.4) }{
        \foreach \Point in {(-0.2,0.25) , (0,0.25),  (0.2,0.25) }{
            \node[darkgreen] at \Point {\textbullet};
        }

	% triangular prism
	\draw[dashed] (-0.75,0) -- (0,0.5) edge (0,2.5) -- (0.75,0);
	\draw[thick] (-0.75,0) rectangle (0.75,2) -- (0,2.5) -- (-0.75,2);
	\end{tikzpicture}}
\caption{Visual representation of the ROD stencil: on the left, a constrained spatial reconstruction is performed at each sub-time node of the ADER correction step; 
on the right, one single space-time constrained polynomial is build by using all the information coming from the space-time predictor. 
{\color{blue}\textbullet} -- quadrature points in the space-time element; 
{\color{darkgreen}\textbullet} -- polynomial coefficients of the local space-time predictor; 
{\color{red}\textbullet} -- sampling points on the real boundary.}
\label{fig:rodSpace_rodSpaceTime}
\end{figure}
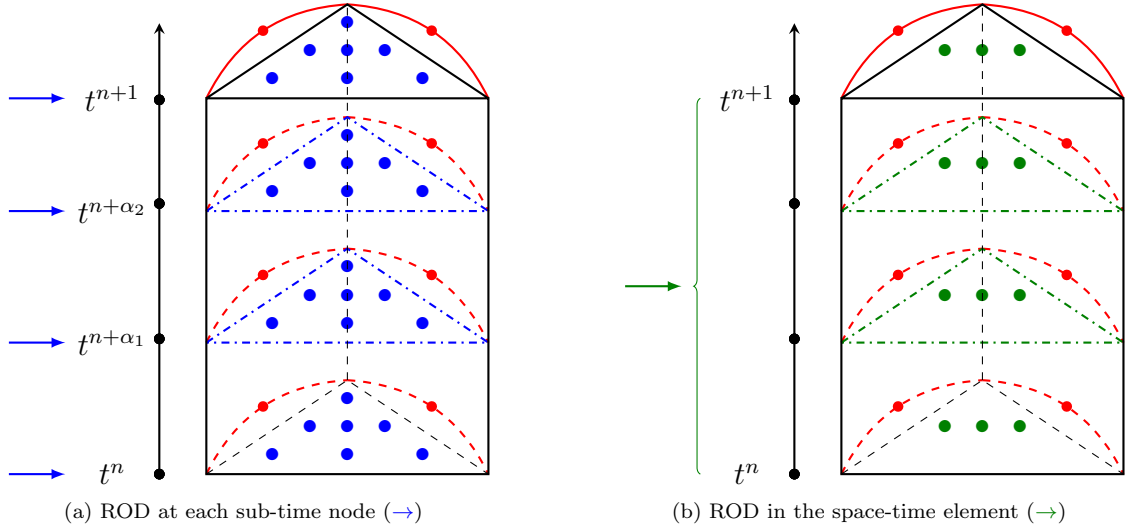

In particular, for all boundary cells a new polynomial in space is defined as,
\begin{equation}
 \varphi(\mathbf{x},\eta,\mathbf{m}) 
= \mathcal{P}_M(\mathbf{x}-\mathbf{m}) \eta 
= \sum_{\alpha=0}^M\sum_{\beta=0}^{M-\alpha}(x_1-m_1)^\alpha(x_2-m_2)^\beta\eta_{(\alpha,\beta)},
\end{equation}
where $\eta$ represents a set of coefficients, and $\mathcal{P}_M(\mathbf{x})$ is a vector of polynomial basis of degree $M$.
In this section, we consider a set of traditional basis functions that differ from those used for the PDE discretization,
\begin{equation}
\begin{aligned}
\setcounter{MaxMatrixCols}{20}
\mathcal{P}_1(x) & = \begin{bmatrix} 1  & x_1 & x_2 \end{bmatrix}, \\
\mathcal{P}_2(x) & = \begin{bmatrix} 1  & x_1 & x_2 & x_1^2 & x_1 x_2 & x_2^2 \end{bmatrix}, \\
\mathcal{P}_3(x) & = \begin{bmatrix} 1  & x_1 & x_2 & x_1^2 & x_1 x_2 & x_2^2 & x_1^3 & x_1^2 x_2 &x_1 x_2^2 & x_2^3\end{bmatrix}.
\end{aligned}
\end{equation}
The point $\mathbf{m}$ is a reference point used for conditioning purposes, and coincides with the barycenter of the element.
It is well-known that a two-dimensional polynomial reconstruction of order $M$ requires $ N =  (M+1)(M+2)/2 $ coefficients. 
\begin{remark}[Comparison with the original approach]
In the ROD initial formulation for finite volume schemes~\cite{costa2018very}, the size of the stencil was chosen $s\approx 1.5 \,N$, and
the corresponding neighboring elements were found using a \textit{searching algorithm}.\\
\end{remark}
When working in the DG framework, a local high order polynomial comes from the discretization, hence there is no need of finding the neighboring elements.
The reconstructed ROD polynomial can easily embed both the information coming from the high order local polynomial and boundary conditions.
It should be noticed that in this case, due to the ROD formulation in space only, the minimization problem must be formulated at each sub-time node of 
the considered time integration scheme, as shown in~\cref{fig:rodSpace_rodSpaceTime}a. 
  
To build a set of sampling points required for the minimization problem, we exploit the local predictor and 
evaluate it in the space quadrature/collocation points $\mathbf{x}_q$ of the sub-time node $t^{n+\alpha_i}$ considered (refer again to~\cref{fig:rodSpace_rodSpaceTime}a).
This gives us a local stencil of sampling points $\bar{\mathcal{U}}_q$ representing fully the high order behavior of the solution at $t^{n+\alpha_i}$.  

Then, the ROD polynomial reconstruction consists in finding the polynomial coefficients $\hat{\eta}$ by minimizing a weighted cost 
function $\hat{\mathcal{J}}$ subject to certain linear constraints  $\hat{\mathcal{G}}$ (the boundary conditions), that reads
\begin{equation}
\begin{aligned}
\hat{\mathcal{J}} (\eta)&= \frac12\sum_{q} \left[w_q \, \mathcal{P}_M(\mathbf{x}_q-\mathbf{m}) \eta - w_q \, \bar{\mathcal{U}}_q \right]^2, \\ 
\text{such that}\quad \hat{\eta}&=\argmin_{\eta} \hat{\mathcal{J}}(\eta)
\quad \text{subject to}\quad \hat{\mathcal{G}}(\eta) = 0.
\end{aligned}
\end{equation}
It should be noticed that the mismatch between the basis functions $\mathcal{P}_M$ used for the ROD reconstruction and those used for the space discretization that define $\bar{\mathcal{U}}$
entails that a least squares problem in terms of quadrature points must be solved.
In the unconstrained case, this can be solved in matrix form as follows,
\begin{equation}
\underbrace{\left[\sum_{q} w_q\mathcal{P}^T_M(\mathbf{x}_q-\mathbf{m})\,
\mathcal{P}_M(\mathbf{x}_q-\mathbf{m})\right]}_{\mathcal{A}_{N\times N}} \, \eta  = 
\underbrace{\left[\sum_{q} w_q\mathcal{P}^T_M(\mathbf{x}_q-\mathbf{m})\, \bar{\mathcal{U}}_q \right]}_{b_{N \times 1}}
\end{equation}
where $w_q:=w(d_q,\delta,\gamma)$ is a {\it weight} that depends on the position of the sampling point $\mathbf{x}_q$.
Usually it is set as the inverse of a distance to a power of $\gamma$:
$$ w(D_q,\delta,\gamma) = \frac1{(\delta D_q)^\gamma + 1} ,$$
where $\delta$ is a sensibility factor and $D_q = |\mathbf{x}^\star-\mathbf{x}_q|$ is the Euclidean distance between point $\mathbf{x}_q$
and the considered quadrature point $\mathbf{x}^\star$.\\ 
%The other parameters are set $\gamma=2$ and $\delta=5|x^\star - m|$. \\

In the constrained case, $\hat\eta$ must fulfill $p$ linear constraints, where $0 < p < N$.
The ROD polynomial is finally obtained by applying these linear constraints in the matrix form $\mathcal{C}\eta=\varrho$,
where $\mathcal{C}$ represents the general Robin boundary conditions $\alpha \varphi + \beta \partial_n \varphi$ imposed on the polynomials $\mathcal{P}_M$
computed on the real boundary $\Gamma=\partial\Omega$, and $\varrho$ is the set of conditions to enforce.
These constraints depend on the nature of the boundary conditions and its geometry through the local normal $n$.  
Following this formulation, we can easily reformulate the linear constraints as follows,
\begin{equation}
\underbrace{\begin{bmatrix} \alpha_{\Gamma_1} \,\mathcal{P}_M(\mathbf{x}_{\Gamma_1}-\mathbf{m}) + \beta_{\Gamma_1} \,\nabla \mathcal{P}_M(\mathbf{x}_{\Gamma_1}-\mathbf{m})\cdot \mathbf{n} \\ \vdots \\ \alpha_{\Gamma_p} \,\mathcal{P}_M(\mathbf{x}_{\Gamma_p}-\mathbf{m}) + \beta_{\Gamma_p} \,\nabla \mathcal{P}_M (\mathbf{x}_{\Gamma_p}-\mathbf{m})\cdot \mathbf{n} \end{bmatrix}}_{\mathcal{C}_{p\times N}} \eta = 
\underbrace{\begin{bmatrix} g_{\Gamma_1} \\ \vdots \\ g_{\Gamma_p} \end{bmatrix}}_{\varrho_{p\times 1}},
\end{equation}
where point $\mathbf{x}_{\Gamma_p}$ is defined as the projection of the quadrature point $\mathbf{x}^\star$ on the real boundary $\partial\Omega$, along the normal $\mathbf{n}$. 
It is straightforward then to retrieve classical Dirichlet ($\beta_{\Gamma_p}=0$) and Neumann ($\alpha_{\Gamma_p}=0$) boundary conditions.\\

As mentioned above this constrained optimization problem can be solved using Lagrange multipliers.
Therefore, we define the functional 
\begin{equation}
\mathcal{L}(\eta,\lambda) = 
\frac12\sum_{q} \left[w_q \, \mathcal{P}_M(\mathbf{x}_q-\mathbf{m}) \eta - w_q \, \bar{\mathcal{U}}_q \right]^2 +
\sum_{k=1}^p \lambda_k \left[ \mathcal{C}\eta - \varrho\right]. 
\end{equation}
The result of this minimization problem $\hat{\eta}$ is the solution of the linear system 
$$\nabla_{\eta,\lambda}\mathcal{L}(\eta,\lambda) =0,$$
 where the differential operator $\nabla_{\eta,\lambda}$ is built by deriving with respect to each
polynomial coefficient and Lagrange multiplier:
\begin{equation}
\left\{
\begin{aligned}
&\nabla_{\eta}    \mathcal{L}(\eta,\lambda) = \mathcal{A}\eta - b + \mathcal{C}^T\lambda=0\\
&\nabla_{\lambda} \mathcal{L}(\eta,\lambda) = \mathcal{C}\eta - \varrho = 0.
\end{aligned}
\right.
\end{equation}
This problem can be solved in matrix form, as follows
\begin{equation}\label{eq:systemLSconstrained}
\left[\begin{array}{c|c} \mathcal{A}  & \mathcal{C}^T \\\hline  \mathcal{C} & 0 \end{array}\right]  
\left[\begin{array}{c} \eta \\\hline \lambda \end{array}\right] = 
\left[\begin{array}{c} b \\\hline \varrho \end{array}\right] 
\qquad \Longrightarrow \qquad
\left[\begin{array}{c} \hat{\eta} \\\hline \hat{\lambda} \end{array}\right] = 
\left[\begin{array}{c|c} \mathcal{A} & \mathcal{C}^T \\\hline \mathcal{C} & 0 \end{array}\right]^{-1} 
\left[\begin{array}{c} b \\\hline \varrho \end{array}\right],
\end{equation}
where $\hat{\eta}$ represents the polynomial coefficients associated to the Lagrange multipliers $\hat{\lambda}$, which are also a
solution of the linear system, that satisfies the boundary conditions.

\subsection{Reconstruction for off-site data methods for space-time schemes}

The original ROD method presented in the previous section is formulated only in space, which entails that the minimization problem
must be solved at each sub-time node $t^{n+\alpha_i}$ of the correction step. 
A similar thing would also happen if dealing with Runge-Kutta methods, meaning that an inversion of the linear system is needed at each 
sub-time node for every boundary element.
In this section, we introduce a novel space-time formulation for the ROD reconstructed polynomial that can easily be integrated 
in the correction step of the ADER method.
Here, we exploit the ADER formulation based on the previous evaluation of the high order space-time predictor~\eqref{eq:predictor0} 
to design a new minimization problem that use the same space-time polynomial basis.
This idea not only allows to write a much simpler optimization problem but also implies the reconstruction of a single space-time ROD polynomial 
(meaning a single inversion of the linear system) used for the boundary conditions of the whole space-time element $K_j\times[t^n,t^{n+1}]$. 

In this case, the fact that the same polynomial basis of the discretization is also used to perform the ROD boundary reconstruction allows us to
write a simplified cost function only with respect to the polynomial coeffiecients, that no longer depends on the evaluation of the solution within the quadrature points, 
\begin{equation}
\begin{aligned}
&\hat{\mathcal{J}} (\eta)= \frac12\left[ \eta  - \hat q \right]^2, \\
\text{such that}\quad \hat{\eta}=&\argmin_{\eta} \hat{\mathcal{J}}(\eta)  \quad \text{subject to}\quad \hat{\mathcal{G}}(\eta) = 0,
\end{aligned}
\end{equation}
where $\hat q$ are the predictor coefficients found with~\cref{eq:predictor0}.

The functional embedding the contratints can now be simplified as 
\begin{equation}
\mathcal{L}(\eta,\lambda) = \frac12\left[ \eta  - \hat q \right]^2 + \sum_{k=1}^p \lambda_k \left[ \mathcal{C}\eta - \varrho \right], 
\end{equation}
from which we can deduce the system minimum as,
\begin{equation}
\begin{cases}
\nabla_{\eta}    \mathcal{L}(\eta,\lambda) = \mathds{1} \eta-\hat{q} + \mathcal{C}^T \lambda=0\\
\nabla_{\lambda} \mathcal{L}(\eta,\lambda) = \mathcal{C}\eta - \varrho = 0.
\end{cases}
\end{equation}
The previous system can be solved once again in matrix form, as follows
\begin{equation}\label{eq:systemLSconstrained2}
\left[\begin{array}{c|c} \mathds{1}  & \mathcal{C}^T \\\hline  \mathcal{C} & 0 \end{array}\right]  
\left[\begin{array}{c} \eta \\\hline \lambda \end{array}\right] = 
\left[\begin{array}{c} \hat{q} \\\hline \varrho \end{array}\right] 
\qquad \Longrightarrow \qquad
\left[\begin{array}{c} \hat{\eta} \\\hline \hat{\lambda} \end{array}\right] = 
\left[\begin{array}{c|c} \mathds{1} & \mathcal{C}^T \\\hline \mathcal{C} & 0 \end{array}\right]^{-1} 
\left[\begin{array}{c} \hat{q}  \\\hline \varrho \end{array}\right].
\end{equation}
By using the space-time polynomial basis we are able to simplify the minimization problem, as now $\mathcal A = \mathds 1$.
\begin{remark}
It should be noticed that the linear system~\eqref{eq:systemLSconstrained2} is slighly larger than the previous one~\eqref{eq:systemLSconstrained}, but it is more diagonally-dominant, and has to be solved only once, which makes it more efficient to invert.
This is also due to the fact that the total number of modes reconstructed is much lower than in the case of several independent reconstructions.
Moreover, by reconstructing a single polynomial for the entire space-time boundary element, the number of required constraints to achieve consistency is halved.
\end{remark}

\subsection{Remarks on the implementation of slip-wall boundary conditions}

In this section we address some implementation features needed to develop the slip-wall boundary condition~\eqref{eq:slipeqv} 
with both the space-only and space-time ROD reconstructions.

\subsubsection{Change of variables}

As mentioned above, for both the original and space-time ROD reconstructions a change of variables is needed to enforce the mixed Dirichlet-Neumann-Robin boundary conditions, meaning that
the ROD polynomial reconstruction must be performed in terms of $\mathcal V(\mathcal U) = (u_n,\,S,\,H,\,u_\tau)$ so that the conditions can be enforced linearly.
When using the original ROD approach, the change of variables is straightforward due to the local evaluation of the sampling points within the 
cell quadrature points $\mathbf{x}_q$. 

However, when working with the space-time formulation only the space-time polynomial coefficients are needed for the minimization problem.
In general such coefficients are defined with respect to the vector of conserved variables, meaning that a projection is needed to retrieve them
for the new vector of variables $\mathcal V$. 

An accurate and efficient way to evaluate them is to perform a space-time $L_2$--projection, 
\begin{equation}\label{eq:proj}
\left(\int_{t^n}^{t^{n+1}} \int_{K_j} \theta_i(\mathbf{x},t) \theta_\ell(\mathbf{x},t)\;\diff{\mathbf{x}}\diff{t}\right)\, \hat{\mathcal V}_\ell  =  \int_{t^n}^{t^{n+1}} \int_{K_j} \theta_i(\mathbf{x},t) \, \mathcal V(\mathcal U_h(\mathbf{x},t)) \;\diff{\mathbf{x}}\diff{t},  
\end{equation}
such that $\mathcal V_h$ can be defined as, 
\begin{equation}\label{eq:slipvariables}
\mathcal V_h(\mathbf{x},t)|_{K_j\times[t^n, t^{n+1}]} = \sum_{\ell}\theta_\ell(\mathbf{x},t)\,\hat{\mathcal V}_\ell.
\end{equation}
\begin{remark}
It should be noticed that, in order to solve problem~\eqref{eq:proj}, $\mathcal V(\mathcal U_h(\mathbf{x},t))$ is evaluated in the space-time
element making necessary the knowledge of the boundary normals and tangents also within the computational domain. 
In our case, since the interface is defined from a level set field, the normal and tangents can be easily computed also far from the boundary.
\end{remark}
\subsubsection{Shock limiting}

For both the ROD formulations, when the \textit{a posteriori} limiter marks a cell on the wall boundary as troubled, 
the boundary conditions must be applied differently because the high order polynomials are discarded, 
and the solution is updated with a second order finite volume scheme.

For simplicity, we adopt a similar flux to the one introduced in~\cite{krivodonova2006high}, which is computed including the information 
related to the real boundary normal.
We take the decomposition $\tilde{\mathbf{n}} = (\tilde{\mathbf{n}} \cdot \mathbf{n} )\mathbf{n} + (\tilde{\mathbf{n}} \cdot \boldsymbol{\tau} )\boldsymbol{\tau}$, where $\tilde{\mathbf{n}}$ represents the normal of the discretized
curved boundary $\partial\Omega_h$ and, $\mathbf{n}$ and $\boldsymbol\tau$ represent the real normal and tangent vectors to $\partial\Omega$, respectively. 
Therefore by writing the flux on the discretized boundary $\mathcal F^{bc}_{\tilde n} = (\tilde{\mathbf{n}} \cdot \mathbf{n} )\mathcal F_n + (\tilde{\mathbf{n}} \cdot \boldsymbol{\tau}  )\mathcal F_\tau$,
and taking $u_n = 0$, we can recast  
\begin{equation}\label{eq:BergerFlux}
\mathcal{F}^{bc}_{\tilde n} = \left(\begin{array}{c} (\tilde{\mathbf{n}} \cdot \boldsymbol{\tau} ) \rho u_\tau \\ (\tilde{\mathbf{n}} \cdot \boldsymbol{\tau} ) \rho \mathbf{u} u_\tau + p \tilde{\mathbf{n}}  \\ (\tilde{\mathbf{n}} \cdot \boldsymbol{\tau} ) \rho u_\tau H  \end{array}\right).
\end{equation}
We remind that~\cref{eq:BergerFlux} is at most third order accurate for general geometries discretized with conformal linear meshes~\cite{costa2018very}, which makes it consistent for troubled cells.

\section{Numerical experiments}
\label{sec:results}

In this section, we present the numerical experiments performed to validate the boundary reconstructions 
presented in~\cref{sec:boundarynumerics}.
In the first example we validate the high order convergence properties of the boundary treatments with 
a smooth manufactured solution characterized by a curved boundary. 
The adaptation of the space-only ROD reconstruction in the ADER-DG framework will be referred to as ROD-DG;
while the novel formulation based on the ADER space-time predictor will be referred to as ROD-ADER-DG.
The numerical results obtained without any correction will be simply referred to as DG.
For the manufactured solution test we present the numerical solutions obtained by imposing Dirichlet-type boundary conditions 
on the curved boundary discretized with both conformal linear elements and fully embedded meshes.
After that, we present several numerical experiments showing the capabilities of the new approach to deal 
with both Dirichlet and slip-wall boundary conditions, also when coupled with an {\it a posteriori} shock limiter. 
Once proven for the first example that the two formulations provide very similar discretization errors, 
for the other test cases we only focus on the results obtained with our new ROD-ADER-DG method.     

\subsection{Manufactured solution: comparisons between the original and space-time reconstruction}\label{test:manuf}

\begin{figure}
\centering
\includegraphics[width=0.45\textwidth]{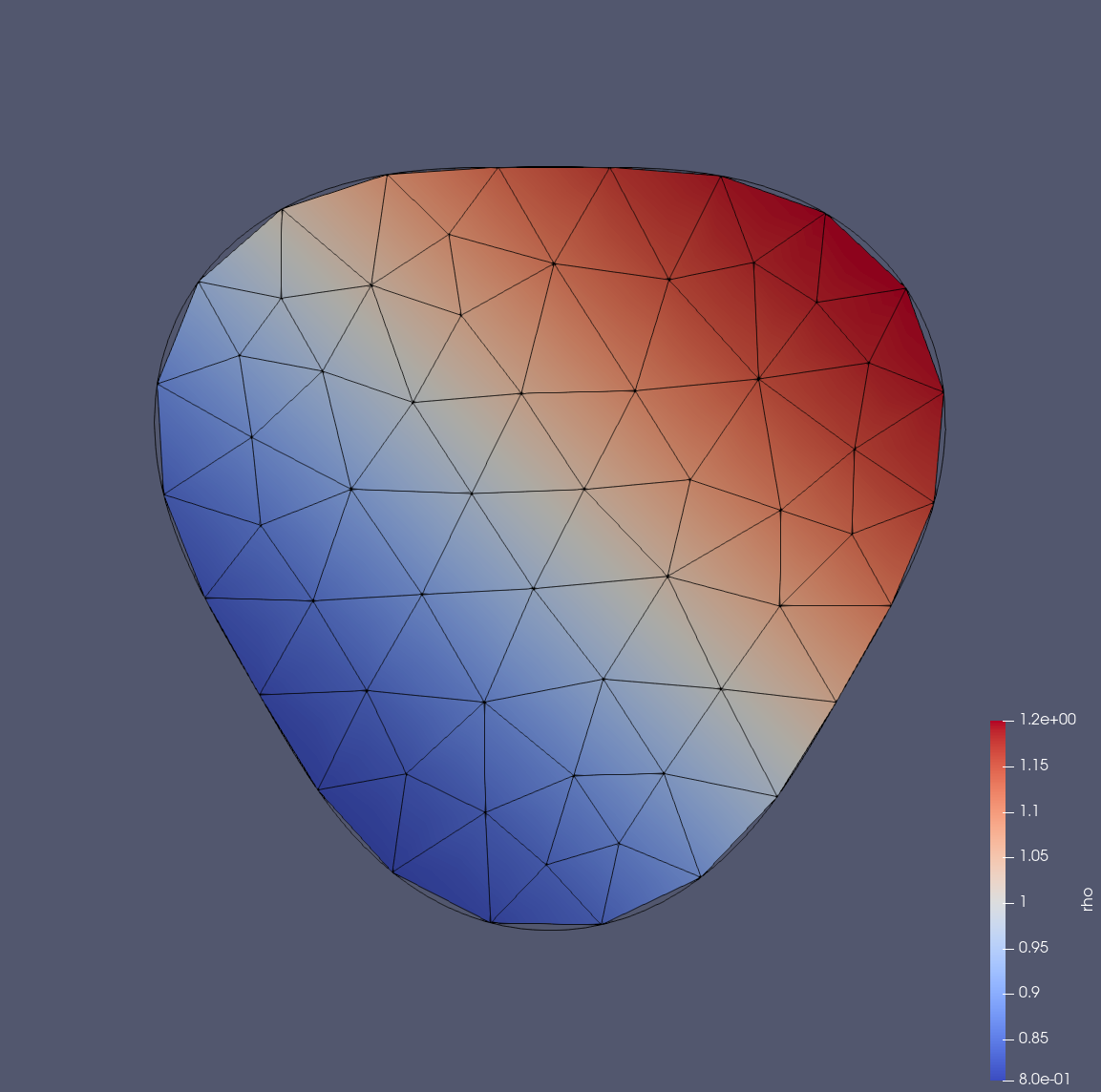}
\caption{Manufactured solution: setup for the test case presented in~\cref{test:manuf}. We plot the coarsest employed conformal mesh, and the initial density profile.}\label{fig:manufconformalmesh}
\end{figure}

The first test case is dedicated to the comparisons between the two ROD formulations described in~\cref{sec:boundarynumerics}.
In particular we focus on a smooth manufactured solution by considering the two-dimensional
inhomogeneous Euler equations~\eqref{eq:euler1} with,
$$\mathcal{S} = (  0.4\cos(x+y),\,0.6\cos(x+y) ,\, 0.6\cos(x+y) ,\, 1.8\cos(x+y))^T.$$
This system has the following exact steady state solution, 
$$\rho = 1+0.2\,\sin(x+y), \quad u = (1,\, 1), \quad p = 1+0.2\,\sin(x+y),$$
which is imposed on the real curved boundary as Dirichlet-type boundary conditions.
A very coarse mesh is generated and then refined by splitting. 
The four grids described in~\cref{tab:manuflinearmesh} are used to study the convergence properties of the formulations.
The boundary conditions are imposed on a complicated curved geometry (see~\cref{fig:manufconformalmesh}), 
that can be described with the following equation written in polar coordinates:
\begin{equation}
\left(\begin{array}{c} x \\ y \end{array}\right) = r(\alpha,\theta) \left(\begin{array}{c} \cos\theta \\ \sin\theta \end{array}\right),\quad \text{ where }\quad r(\alpha,\theta)=r_0\left(1+\frac{1}{10}\sin(\alpha\theta)\right),
\end{equation}
where $r_0=1$ and $\alpha=3$.

\begin{table}
        \caption{Manufactured solution: characteristics of the conformal linear meshes for the test case presented in~\cref{test:manuf}.}
        \label{tab:manuflinearmesh}
        \footnotesize
        \centering
        \begin{tabular}{cccc} \hline %\hline
                Grid level &Nodes  &Triangles  &$h$ \\[0.5mm]
                \hline
                0 & 61     &  98    & 1.54E-1   \\
                1 & 219    &  392   & 7.72E-2   \\
                2 & 829    & 1,568  & 3.86E-2   \\
                3 & 3,225  & 6,272  & 1.93E-2   \\
                \hline %\hline
        \end{tabular}
\end{table}

Since the standard Dirichlet boundary condition is enforced onto a curved boundary, discretized with a polygonal mesh, 
for the DG solution we expect to have at most second order of accuracy no matter the order of the polynomial used.
This is clear from~\cref{tab:manufnocorrection} where the $L_2$ error norm is shown for the different experiments.

\begin{table}
\caption{Manufactured solution: convergence analysis for the test case presented in~\cref{test:manuf}. Numerical results obtained with Dirichlet-type boundary conditions without any correction on the conformal linear meshes described in~\cref{tab:manuflinearmesh}. One can notice that by imposing naively the boundary conditions, the method is at most second order accurate.}\label{tab:manufnocorrection}
\footnotesize
\centering
\begin{tabular}{ccccccccc}
        \hline
        &\multicolumn{2}{c}{$\rho$} &\multicolumn{2}{c}{$\rho u$} &\multicolumn{2}{c}{$\rho v$}  &\multicolumn{2}{c}{$\rho E$}\\[0.5mm]
        \cline{2-9}
        Grid level & $L_2$        & $\tilde{n}$ & $L_2$        & $\tilde{n}$ & $L_2$        & $\tilde{n}$ & $L_2$      & $\tilde{n}$ \\[0.5mm]\hline
        &\multicolumn{8}{c}{DG-$\mathcal{P}_1$}\\[0.5mm]
        0          &  1.7293E-3   &   --        &  2.3514E-3   &  --         &  2.2344E-3   &   --        &  6.6466E-3 &   --    \\
        1          &  4.4790E-4   &    1.95     &  6.0925E-4   &    1.95     &  5.7869E-4   &    1.95     &  1.7138E-3 &  1.96   \\
        2          &  1.1276E-4   &    1.99     &  1.5358E-4   &    1.99     &  1.4579E-4   &    1.99     &  4.3195E-4 &  1.99   \\
        3          &  2.8241E-5   &    2.00     &  3.8473E-5   &    2.00     &  3.6512E-5   &    2.00     &  1.0820E-4 &  2.00   \\
        &\multicolumn{8}{c}{DG-$\mathcal{P}_2$}\\[0.5mm]
        0          &  1.6583E-3   &   --        &  2.3344E-3   &  --         &  2.2040E-3   &   --        &  6.4258E-3 &   --    \\
        1          &  4.2195E-4   &   1.97      &  5.9156E-4   &    1.98     &  5.5908E-4   &    1.98     &  1.6308E-3 &  1.98   \\
        2          &  1.0522E-4   &   2.00      &  1.4777E-4   &    2.00     &  1.3975E-4   &    2.00     &  4.0791E-4 &  2.00   \\
        3          &  2.6187E-5   &   2.01      &  3.6855E-5   &    2.00     &  3.4867E-5   &    2.00     &  1.0170E-4 &  2.00   \\
        &\multicolumn{8}{c}{DG-$\mathcal{P}_3$}\\[0.5mm]
        0          &  1.6817E-3   &   --        &  2.3809E-3   &  --         &  2.2562E-3   &   --        &  6.5431E-3 &   --    \\
        1          &  4.3095E-4   &   1.96      &  6.0642E-4   &   1.97      &  5.7430E-4   &   1.97      &  1.6640E-3 &  1.98   \\
        2          &  1.0768E-4   &   2.00      &  1.5092E-4   &   2.01      &  1.4280E-4   &   2.01      &  4.1580E-4 &  2.00   \\
        3          &  2.6792E-5   &   2.01      &  3.7467E-5   &   2.01      &  3.5428E-5   &   2.01      &  1.0345E-4 &  2.01   \\
        \hline\\[1pt]
\end{tabular}
\end{table}

\begin{figure}
\centering
\subfloat[Mesh $5\times 5$]  {\includegraphics[width=0.30\textwidth]{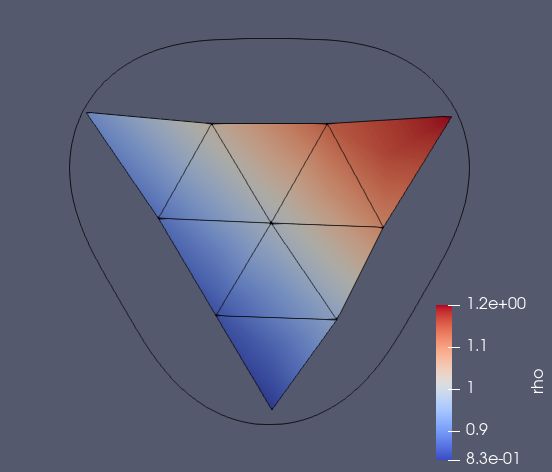}}\quad
\subfloat[Mesh $10\times 10$]{\includegraphics[width=0.30\textwidth]{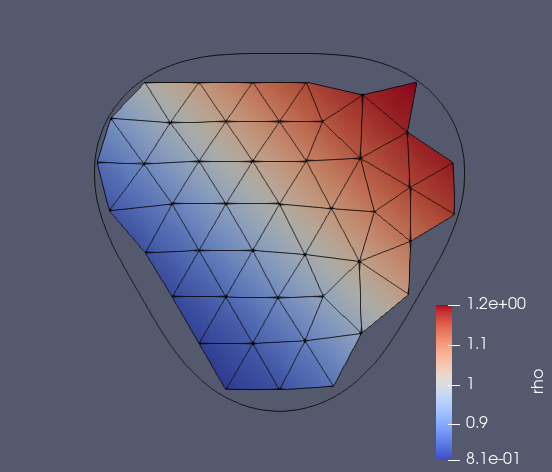}}\quad
\subfloat[Mesh $15\times 15$]{\includegraphics[width=0.30\textwidth]{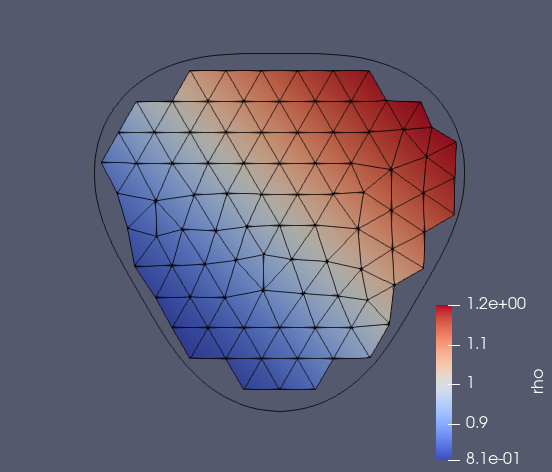}}\qquad
\subfloat[Mesh $20\times 20$]{\includegraphics[width=0.30\textwidth]{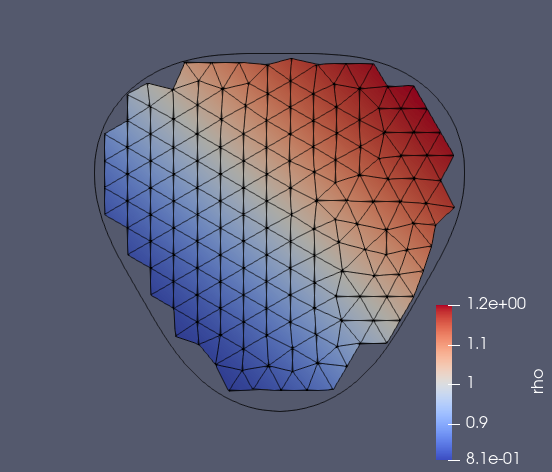}}\qquad
\subfloat[Mesh $35\times 35$]{\includegraphics[width=0.30\textwidth]{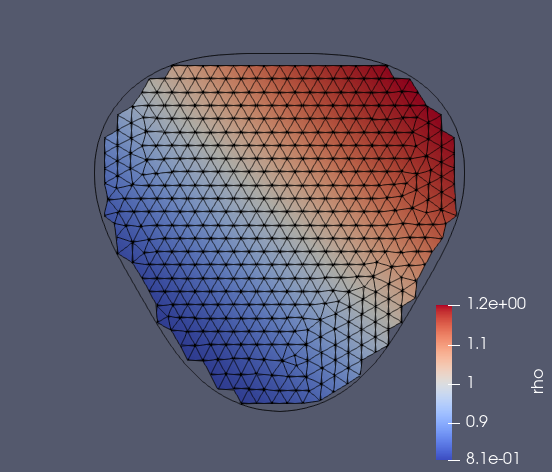}}\qquad
\caption{Manufactured solution: setup for the test case presented in~\cref{test:manuf}. We plot the employed embedded linear meshes on a square $[-1.5,1.5]\times[-1.5,1.5]$, and the initial density profile.}\label{fig:manufembeddedmesh}
\end{figure}

We then perform the same analysis with both the original and space-time ROD formulations to prove numerically their high order 
convergence properties. 
We present the numerical results in~\cref{tab:manuforiginalRODconformal,tab:manuforiginalRODembedded,tab:manufspacetimeRODconformal,tab:manufspacetimeRODembedded} 
obtained with both formulations on conformal linear meshes and fully embedded meshes. 
In all cases, the expected convergence of $M+1$ is achieved. 
Moreover, as mentioned above, the new space-time formulation is able to achieve similar results while being less computationally expensive. 
In particular, we experience a reduction of $\sim 10\%$ for $\mathcal P_3$ polynomials, which is remarkable given that the ROD reconstruction only influences the boundary cells.
It should be noticed that the convergence trends obtained by the method on fully embedded meshes 
(described in~\cref{tab:manufembedmesh}) are typical of such methods~\cite{Scovazzi1}, simply because by refining 
the mesh the discretized  boundary gets closer to the real one depending on how it cuts the background mesh 
(see~\cref{fig:manufembeddedmesh}).  

\begin{table}
	\caption{Manufactured solution: characteristics of the fully embedded meshes on a square $[-1.5,1.5]\times[-1.5,1.5]$ for the test case presented in~\cref{test:manuf}.}
	\label{tab:manufembedmesh}
	\footnotesize
	\centering
	\begin{tabular}{cccc} \hline %\hline
		Grid level  &$h$ \\[0.5mm]
		\hline
		$ 5 \times  5$  & 3.79E-1   \\
		$10 \times 10$  & 2.02E-1   \\
		$15 \times 15$  & 1.33E-1   \\
		$20 \times 20$  & 9.93E-2   \\
		$35 \times 35$  & 6.02E-2   \\ 
		\hline %\hline
	\end{tabular}
\end{table}

\begin{table}
        \caption{Manufactured solution: convergence analysis for the test case presented in~\cref{test:manuf}. Numerical results obtained with the ROD-DG reconstruction on the conformal linear meshes described in~\cref{tab:manuflinearmesh}. One can notice that the high order precision is achieved for all polynomials.}%\\}
        \label{tab:manuforiginalRODconformal}
        \footnotesize
        \centering
        \begin{tabular}{ccccccccc}
                \hline%\\[1pt]
                %\multicolumn{9}{c}{Convergence analysis with the ROD-DG reconstruction on conformal linear meshes}\\
                %\hline
                &\multicolumn{2}{c}{$\rho$} &\multicolumn{2}{c}{$\rho u$} &\multicolumn{2}{c}{$\rho v$}  &\multicolumn{2}{c}{$\rho E$}\\[0.5mm]
                \cline{2-9}
                Grid level & $L_2$        & $\tilde{n}$ & $L_2$        & $\tilde{n}$ & $L_2$        & $\tilde{n}$ & $L_2$      & $\tilde{n}$ \\[0.5mm]\hline
                &\multicolumn{8}{c}{ROD-DG-$\mathcal{P}_1$}\\[0.5mm]
                0              & 9.0353E-4 &   --        & 9.3786E-4 &  --         & 9.3503E-4 &   --        & 2.8893E-3  &   --     \\
                1              & 2.0990E-4 &   2.11      & 2.0792E-4 &   2.17      & 2.0892E-4 &    2.16     & 6.6542E-4  &   2.12   \\
                2              & 5.0528E-5 &   2.05      & 4.9278E-5 &   2.08      & 4.9736E-5 &    2.07     & 1.6082E-4  &   2.05   \\
                3              & 1.2438E-5 &   2.02      & 1.2032E-5 &   2.03      & 1.2173E-5 &    2.03     & 3.9724E-5  &   2.02   \\
                &\multicolumn{8}{c}{ROD-DG-$\mathcal{P}_2$}\\[0.5mm]
                0              & 5.8327E-5 &   --        & 4.0652E-5 &  --         & 4.5945E-5 &   --        & 1.7771E-4  &   --     \\
                1              & 9.5321E-6 &    2.61     & 5.1266E-6 &    2.99     &  6.2677E-6&   2.87      & 2.8402E-5  &   2.65   \\
                2              & 1.4695E-6 &    2.70     & 6.9200E-7 &    2.89     & 8.6795E-7 &   2.85      & 4.3499E-6  &   2.71   \\
                3              & 2.0870E-7 &    2.82     & 9.7264E-8 &    2.83     & 1.1487E-7 &   2.92      & 6.1764E-7  &   2.82   \\
                &\multicolumn{8}{c}{ROD-DG-$\mathcal{P}_3$}\\[0.5mm]
                0              & 1.0030E-6 &   --        & 1.3555E-6 &  --         & 1.4008E-6 &   --        & 3.4352E-6  &   --     \\
                1              & 4.7755E-8 &   4.39      & 5.4425E-8 &   4.64      & 5.5114E-8 &   4.67      & 1.5203E-7  &   4.50   \\
                2              & 2.5871E-9 &   4.21      & 2.7218E-9 &   4.32      & 2.7433E-9 &   4.33      & 8.1616E-9  &   4.22   \\
                3              & 1.5035E-10&   4.10      & 1.5323E-10&   4.15      & 1.5474E-10&   4.15      & 4.7689E-10 &   4.10   \\
                \hline
        \end{tabular}
\end{table}

\begin{table}
        \caption{Manufactured solution: convergence analysis for the test case presented in~\cref{test:manuf}. Numerical results obtained with the ROD-DG reconstruction on the fully embedded meshes presented in~\cref{fig:manufembeddedmesh}. One can notice that the high order precision is achieved for all polynomials.}%\\}
        \label{tab:manuforiginalRODembedded}
        \footnotesize
        \centering
        \begin{tabular}{ccccccccc}
                \hline%\\[1pt]
                &\multicolumn{2}{c}{$\rho$} &\multicolumn{2}{c}{$\rho u$} &\multicolumn{2}{c}{$\rho v$}  &\multicolumn{2}{c}{$\rho E$}\\[0.5mm]
                \cline{2-9}
                Grid level & $L_2$        & $\tilde{n}$ & $L_2$        & $\tilde{n}$ & $L_2$        & $\tilde{n}$ & $L_2$      & $\tilde{n}$ \\[0.5mm]\hline
                &\multicolumn{8}{c}{ROD-DG-$\mathcal{P}_1$}\\[0.5mm]
                $ 5 \times  5$ &  0.1840E-01  &    --   &  0.2589E-01  &    --   &  0.2369E-01  &    --   &  0.7320E-01  &    --  \\
                $10 \times 10$ &  0.4015E-02  &   2.26  &  0.5704E-02  &   2.24  &  0.5326E-02  &   2.21  &  0.1566E-01  &   2.28 \\
                $15 \times 15$ &  0.2141E-02  &   1.64  &  0.3118E-02  &   1.57  &  0.3067E-02  &   1.44  &  0.8588E-02  &   1.56 \\
                $20 \times 20$ &  0.1732E-02  &   0.74  &  0.2453E-02  &   0.83  &  0.2602E-02  &   0.57  &  0.7137E-02  &   0.64 \\
                $35 \times 35$ &  0.4966E-03  &   2.26  &  0.7377E-03  &   2.17  &  0.7078E-03  &   2.36  &  0.1996E-02  &   2.30 \\
                &\multicolumn{8}{c}{ROD-DG-$\mathcal{P}_2$}\\[0.5mm]
                $ 5 \times  5$ &  0.4464E-02  &    --   &  0.7350E-02  &    --   &  0.6120E-02  &    --   &  0.1918E-01  &    --  \\
                $10 \times 10$ &  0.2611E-03  &   4.21  &  0.3919E-03  &   4.34  &  0.3617E-03  &   4.19  &  0.1029E-02  &   4.33 \\
                $15 \times 15$ &  0.9730E-04  &   2.57  &  0.1514E-03  &   2.48  &  0.1499E-03  &   2.29  &  0.3833E-03  &   2.57 \\
                $20 \times 20$ &  0.5464E-04  &   2.00  &  0.7456E-04  &   2.46  &  0.8333E-04  &   2.04  &  0.2213E-03  &   1.90 \\
                $35 \times 35$ &  0.8168E-05  &   3.44  &  0.1196E-04  &   3.31  &  0.1057E-04  &   3.73  &  0.3169E-04  &   3.52 \\
                &\multicolumn{8}{c}{ROD-DG-$\mathcal{P}_3$}\\[0.5mm]
                $ 5 \times  5$ &  0.5828E-03  &    --   &  0.1035E-02  &    --   &  0.9056E-03  &    --   &  0.2743E-02  &    --  \\
                $10 \times 10$ &  0.3038E-04  &   4.38  &  0.4969E-04  &   4.50  &  0.4695E-04  &   4.38  &  0.1252E-03  &   4.57 \\
                $15 \times 15$ &  0.7848E-05  &   3.52  &  0.1435E-04  &   3.23  &  0.1368E-04  &   3.21  &  0.3360E-04  &   3.43 \\
                $20 \times 20$ &  0.4749E-05  &   1.74  &  0.7741E-05  &   2.14  &  0.7982E-05  &   1.87  &  0.2019E-04  &   1.77 \\
                $35 \times 35$ &  0.4305E-06  &   4.34  &  0.7545E-06  &   4.21  &  0.7202E-06  &   4.35  &  0.1822E-05  &   4.35 \\
                \hline
        \end{tabular}
\end{table}

\begin{table}
        \caption{Manufactured solution: convergence analysis for the test case presented in~\cref{test:manuf}. Numerical results obtained with the ROD-ADER-DG reconstruction on the conformal linear meshes described in~\cref{tab:manuflinearmesh}. One can notice that the high order precision is achieved for all polynomials.}%\\}
        \label{tab:manufspacetimeRODconformal}
        \footnotesize
        \centering
        \begin{tabular}{ccccccccc}
                \hline%\\[1pt]
                %\multicolumn{9}{c}{Convergence analysis with the ROD-ADER-DG reconstruction on conformal linear meshes}\\
                %\hline
                &\multicolumn{2}{c}{$\rho$} &\multicolumn{2}{c}{$\rho u$} &\multicolumn{2}{c}{$\rho v$}  &\multicolumn{2}{c}{$\rho E$}\\[0.5mm]
                \cline{2-9}
                Grid level & $L_2$    & $\tilde{n}$ & $L_2$   & $\tilde{n}$ & $L_2$   & $\tilde{n}$ & $L_2$ & $\tilde{n}$ \\[0.5mm]\hline
                &\multicolumn{8}{c}{ROD-ADER-DG-$\mathcal{P}_1$}\\[0.5mm]
                0      &  8.3292E-4   &    --    &  8.2968E-4   &    --     &  8.2849E-4   &   --     &  2.6472E-3   &   --       \\
                1      &  2.0273E-4   &   2.04   &  1.9689E-4   &   2.07    &  1.9822E-4   &   2.06   &  6.4251E-4   &   2.04     \\
                2      &  4.9613E-5   &   2.03   &  4.7793E-5   &   2.04    &  4.8286E-5   &   2.04   &  1.5771E-4   &   2.03     \\
                3      &  1.2295E-5   &   2.01   &  1.1778E-5   &   2.02    &  1.1921E-5   &   2.02   &  3.9159E-5   &   2.01     \\
                &\multicolumn{8}{c}{ROD-ADER-DG-$\mathcal{P}_2$}\\[0.5mm]
                0      &  5.7478E-5   &    --    &  3.9499E-5   &    --     &  4.4647E-5   &   --     &  1.7488E-4   &   --       \\
                1      &  9.4589E-6   &   2.60   &  5.0413E-6   &   2.97    &  6.1664E-6   &   2.86   &  2.8118E-5   &   2.64     \\
                2      &  1.4611E-6   &   2.70   &  6.8355E-7   &   2.88    &  8.5622E-7   &   2.85   &  4.3155E-6   &   2.70     \\
                3      &  2.0771E-7   &   2.81   &  9.6318E-8   &   2.83    &  1.1345E-7   &   2.92   &  6.1345E-7   &   2.81     \\
                &\multicolumn{8}{c}{ROD-ADER-DG-$\mathcal{P}_3$}\\[0.5mm]
                0      &  8.1215E-7   &    --    &  1.0517E-6   &    --     &  1.0822E-6   &   --     &  2.8396E-6   &   --       \\
                1      &  4.3856E-8   &   4.21   &  4.9314E-8   &   4.41    &  5.0113E-8   &   4.43   &  1.4269E-7   &  4.31      \\
                2      &  2.4815E-9   &   4.14   &  2.5941E-9   &   4.25    &  2.6192E-9   &   4.25   &  7.9180E-9   &  4.17      \\
                3      &  1.4723E-10  &   4.07   &  1.4964E-10  &   4.12    &  1.5109E-10  &   4.11   &  4.7004E-10  &  4.07    \\
                \hline
        \end{tabular}
\end{table}

\begin{table}
        \caption{Manufactured solution: convergence analysis for the test case presented in~\cref{test:manuf}. Numerical results obtained with the ROD-ADER-DG reconstruction on the fully embedded meshes presented in~\cref{fig:manufembeddedmesh}. One can notice that the high order precision is achieved for all polynomials.}%\\}
        \label{tab:manufspacetimeRODembedded}
        \footnotesize
        \centering
        \begin{tabular}{ccccccccc}
                \hline%\\[1pt]
                %\multicolumn{9}{c}{Convergence analysis with the ROD reconstruction on embedded meshes on a square $[-1.5,1.5]\times[-1.5,1.5]$}\\
                %\hline
                &\multicolumn{2}{c}{$\rho$} &\multicolumn{2}{c}{$\rho u$} &\multicolumn{2}{c}{$\rho v$}  &\multicolumn{2}{c}{$\rho E$}\\[0.5mm]
                \cline{2-9}
                Grid level & $L_2$        & $\tilde{n}$ & $L_2$        & $\tilde{n}$ & $L_2$        & $\tilde{n}$ & $L_2$      & $\tilde{n}$ \\[0.5mm]\hline
                &\multicolumn{8}{c}{ROD-ADER-DG-$\mathcal{P}_1$}\\[0.5mm]
                $ 5 \times  5$ & 2.1244E-2 &  --    & 2.8921E-2 &   --    & 2.6676E-2 &   --    & 8.2175E-2 &   --      \\ 
                $10 \times 10$ & 4.7323E-3 & 2.39   & 6.9117E-3 &  2.28   & 6.4829E-3 &  2.25   & 1.9135E-2 &  2.32     \\
                $15 \times 15$ & 2.5373E-3 & 1.51   & 3.7614E-3 &  1.47   & 3.7089E-3 &  1.35   & 1.0495E-2 &  1.45     \\
                $20 \times 20$ & 1.9806E-3 & 0.83   & 2.8678E-3 &  0.91   & 3.0124E-3 &  0.70   & 8.3544E-3 &  0.76     \\
                $35 \times 35$ & 5.7732E-4 & 2.46   & 8.7147E-4 &  2.38   & 8.3663E-4 &  2.26   & 2.3870E-3 &  2.50     \\
                &\multicolumn{8}{c}{ROD-ADER-DG-$\mathcal{P}_2$}\\[0.5mm]
                $ 5 \times  5$ & 6.4776E-3 &  --    & 1.0806E-2 &  --     & 8.9916E-3 &   --    & 2.8161E-2 &   --      \\ 
                $10 \times 10$ & 4.3917E-4 &  4.29  & 6.7778E-4 &   4.41  & 6.4805E-4 &  4.18   & 1.8388E-3 &  4.34     \\
                $15 \times 15$ & 1.6884E-4 &  2.31  & 2.6535E-4 &   2.26  & 2.6897E-4 &  2.12   & 7.0304E-4 &  2.32     \\
                $20 \times 20$ & 8.9175E-5 &  2.15  & 1.2627E-4 &   2.50  & 1.3957E-4 &  2.20   & 3.7224E-4 &  2.14     \\
                $35 \times 35$ & 1.3315E-5 &  3.79  & 2.0681E-5 &   3.61  & 1.7738E-5 &  4.12   & 5.4102E-5 &  3.85     \\
                &\multicolumn{8}{c}{ROD-ADER-DG-$\mathcal{P}_3$}\\[0.5mm]
                $ 5 \times  5$ & 7.4636E-4 &  --    & 1.3662E-3 &  --     & 1.1646E-3 &   --    & 3.6090E-3 &   --      \\ 
                $10 \times 10$ & 7.0291E-5 &  3.76  & 1.0730E-4 &  4.05   & 1.0734E-4 &  3.80   & 2.9700E-4 &  3.98     \\
                $15 \times 15$ & 1.2500E-5 &  4.17  & 2.1671E-5 &  3.86   & 2.2876E-5 &  3.73   & 5.7446E-5 &  3.97     \\
                $20 \times 20$ & 7.9376E-6 &  1.53  & 1.3160E-5 &  1.68   & 1.2447E-5 &  2.05   & 3.2858E-5 &  1.88     \\
                $35 \times 35$ & 8.4487E-7 &  4.47  & 1.4199E-6 &  4.45   & 1.3451E-6 &  4.44   & 3.5280E-6 &  4.46     \\
                \hline
        \end{tabular}
\end{table}

\subsection{Isentropic supersonic vortex: validation of Dirichlet-type and slip-wall boundary conditions}\label{test:vort}

We now consider an isentropic supersonic flow between two concentric circular arcs of radii $r_i = 1$ and
$r_o = 1.384$. The exact density in terms of radius $r$ is given by
\begin{equation}\label{eq:vortex2D rho}
\rho = \rho_i \left( 1 + \frac{\gamma+1}{2}M_i^2  \left(1- \left(\frac{r_i}{r}\right)^2 \right) \right)^{\frac{1}{\gamma-1}}.
\end{equation}
The velocity and pressure are given by
\begin{equation}\label{eq:vortex2D u p}
\|u\|=\frac{c_i M_i}{r}, \quad p = \frac{\rho^\gamma}{\gamma},
\end{equation}
where $c_i$ is the speed of sound on the inner circle. The Mach number on the inner circle $M_i$ is set to 2.25 and
the density $\rho_i$ to 1.

\begin{figure}
\centering
\includegraphics[width=0.45\textwidth]{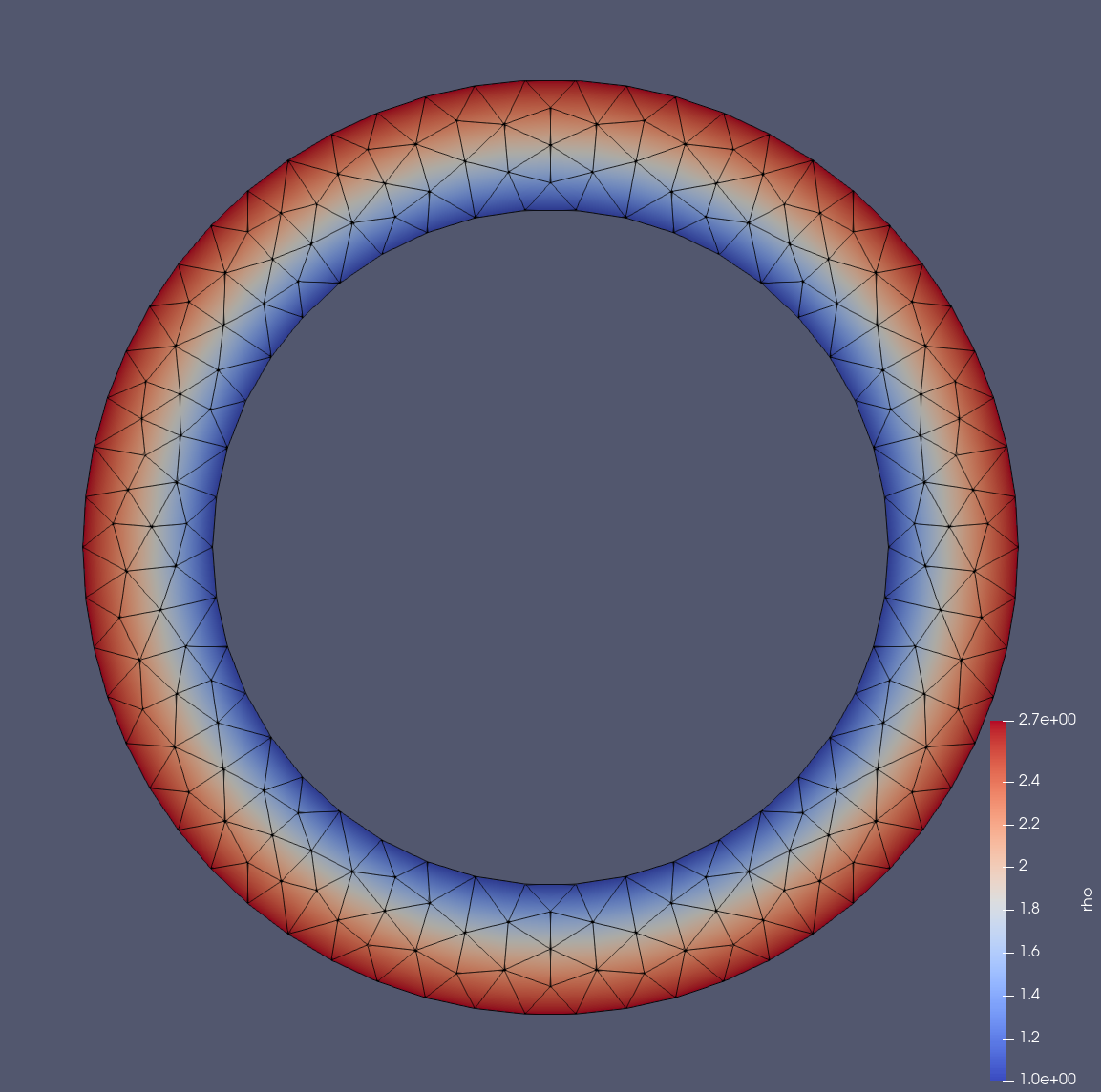}
\caption{Isentropic supersonic vortex: setup for the test case presented in~\cref{test:vort}. We plot a coarse mesh and the initial density profile.}\label{fig:vortexconformalmesh}
\end{figure}

In this section we present only the results obtained with the novel formulation, which has already shown better performances, for different boundary conditions. 
Here we present the results obtained on four refined conformal triangulations. 
In~\cref{fig:vortexconformalmesh} the computational domain and a coarse triangulation are shown.

\begin{table}
        \caption{Isentropic supersonic vortex: characteristics of the conformal linear meshes for the test case presented in~\cref{test:vort}.} 
        \label{tab:meshconfvortex}
        \footnotesize
        \centering
        \begin{tabular}{cccc} \hline %\hline
                Grid level &Nodes  &Triangles & h \\[0.5mm]
                \hline
                0 & 71     & 92      &  1.50E-1 \\
                1 & 238    & 376     &  7.50E-2 \\
                2 & 852    & 1,504   &  3.75E-2 \\
                3 & 3,208  & 6,016   &  1.87E-2 \\
                %3 & 12,432 & 24,064 &          \\
                \hline %\hline
        \end{tabular}
\end{table}

Since for this test case we have an analytical solution, we present in~\cref{tab:vortexconformalnorod} the results 
of the experiment performed by enforcing Dirichlet-type conditions on the internal boundary, 
and classical slip-wall conditions on the external boundary without any correction. 
As expected, since both boundaries are discretized with linear meshes the convergence trends cannot be higher than 2.
Moreover, it is well-known that the reflecting boundary conditions give rise to severe spurious errors when increasing the 
polynomial order~\cite{bassi1997high,krivodonova2006high}. For this reason, on the same mesh, increased errors are found when increasing the order of the discretization.
Similar trends are also shown in other works in the literature~\cite{krivodonova2006high,ciallella2023shifted}. 
Notice that due to important spurious disturbances, the simulation performed with $\mathcal P_3$ on the coarsest mesh does not provide any
result, i.e.\ the simulation accidentally stops before the final time is reached.
We validate again the ROD boundary treatment on this experiment, achieving the expected discretization errors and order-of-accuracy 
(see~\cref{tab:vortexconformalrod}). 
Notice that when enforcing slip-wall conditions, the estimated order-of-accuracy is $M$ rather than $M+1$.
This is due to the fact that the complex slip-wall treatment needs Neumann and Robin boundary conditions that are imposed using 
the derivative of the polynomials, hence losing one order of accuracy.
This is also typical of other similar approaches~\cite{costa2019very,visbech2023spectral}.
Although the boundary treatment loses one order-of-accuracy, all discretization errors found in our experiments are several orders of magnitude smaller than those
obtained with no corrections, which are at most second order. Moreover  increasing the polynomial order always allows a reduction in error of at least one order of magnitude, even when the increase of rate of convergence is mild.  
A trivial way to overcome this issue consists in constructing derivative polynomials of degree M using some informations from neighboring cells. 
This could be done only for those cells on the Neumann/Robin boundary, as shown in~\cite{costa2019very}. 
However, this is out of the scope of this work therefore we will not discuss this case here.

\begin{table}
\caption{Isentropic supersonic vortex: convergence analysis for the test case presented in~\cref{test:vort}. Numerical results obtained with Dirichlet for the internal boundary and slip-wall for the external boundary without any correction on the conformal linear meshes described in~\cref{tab:meshconfvortex}. One can notice that, due to the boundary approximation, the classical slip-wall condition provide convergence trends lower than 2.}\label{tab:vortexconformalnorod}
\footnotesize
\centering
\begin{tabular}{ccccccccc}
        \hline
        &\multicolumn{2}{c}{$\rho$} &\multicolumn{2}{c}{$\rho u$} &\multicolumn{2}{c}{$\rho v$}  &\multicolumn{2}{c}{$\rho E$}\\[0.5mm]
        \cline{2-9}
        Grid level & $L_2$        & $\tilde{n}$ & $L_2$        & $\tilde{n}$ & $L_2$        & $\tilde{n}$ & $L_2$   & $\tilde{n}$ \\[0.5mm]\hline
        &\multicolumn{8}{c}{DG-$\mathcal{P}_1$}\\[0.5mm]
        0     & 1.0757E-1  &    --       &  2.1620E-1 &    --       &  2.1855E-1  &   --        & 4.9486E-1  &   --       \\  
        1     & 3.9747E-2  &    1.44     &  8.4226E-2 &   1.36      &  8.4449E-2  &   1.37      & 1.8268E-1  &   1.44     \\  
        2     & 1.3584E-2  &    1.55     &  2.9191E-2 &   1.53      &  2.9274E-2  &   1.53      & 6.2762E-2  &   1.54     \\  
        3     & 4.8112E-3  &    1.50     &  1.0417E-2 &   1.49      &  1.0448E-2  &   1.49      & 2.2273E-2  &   1.50     \\  
        &\multicolumn{8}{c}{DG-$\mathcal{P}_2$}\\[0.5mm]
        0     & 1.5596E-1  &    --       &  3.6487E-1 &    --       &  3.6579E-1  &   --        & 7.4800E-1  &   --       \\  
        1     & 6.0207E-2  &    1.37     &  1.3306E-1 &   1.46      &  1.3325E-1  &   1.46      & 2.8177E-1  &   1.41     \\  
        2     & 2.1706E-2  &    1.47     &  4.7605E-2 &   1.48      &  4.7647E-2  &   1.48      & 1.0117E-1  &   1.48     \\  
        3     & 7.8804E-3  &    1.46     &  1.6924E-2 &   1.49      &  1.6933E-2  &   1.49      & 3.6319E-2  &   1.48     \\  
        &\multicolumn{8}{c}{DG-$\mathcal{P}_3$}\\[0.5mm]
        0     &   --       &    --       &   --       &    --       &   --        &   --        &  --        &   --       \\  
        1     & 1.4987E-1  &    --       &  1.5568E-1 &    --       &  1.5817E-1  &   --        & 5.5100E-1  &   --       \\  
        2     & 6.9222E-2  &   1.11      &  6.5801E-2 &    1.24     &  6.6299E-2  &   1.25      & 2.5165E-1  &   1.13     \\  
        3     & 3.1405E-2  &   1.14      &  2.8275E-2 &    1.22     &  2.8553E-2  &   1.22      & 1.1342E-1 &    1.15    \\  
        \hline\\[1pt]
\end{tabular}
\end{table}

\begin{table}
\caption{Isentropic supersonic vortex: convergence analysis for the test case presented in~\cref{test:vort}. Numerical results obtained with Dirichlet for the internal boundary and slip-wall for the external boundary with the ROD-ADER-DG reconstruction on the conformal linear meshes described in~\cref{tab:meshconfvortex}. One can notice that the expected sub-optimal precision given by Neumann conditions is achieved.}\label{tab:vortexconformalrod}
\footnotesize
\centering
\begin{tabular}{ccccccccc}
        \hline
        &\multicolumn{2}{c}{$\rho$} &\multicolumn{2}{c}{$\rho u$} &\multicolumn{2}{c}{$\rho v$}  &\multicolumn{2}{c}{$\rho E$}\\[0.5mm]
        \cline{2-9}
        Grid level & $L_2$        & $\tilde{n}$ & $L_2$        & $\tilde{n}$ & $L_2$        & $\tilde{n}$ & $L_2$      & $\tilde{n}$ \\[0.5mm]\hline
        &\multicolumn{8}{c}{ROD-ADER-DG-$\mathcal{P}_1$}\\[0.5mm]
        0     &  1.2539E-1   &    --       &  1.6824E-1   &    --       &  1.6893E-1   &   --        &  5.5124E-1   &   --       \\  
        1     &  3.7048E-2   &    1.76     &  5.4835E-2   &    1.62     &  5.4921E-2   &   1.62      &  1.7740E-1   &   1.64     \\  
        2     &  9.5025E-3   &    1.96     &  1.3583E-2   &    2.01     &  1.3596E-2   &   2.01      &  4.5769E-2   &   1.95     \\  
        3     &  2.2951E-3   &    2.05     &  3.1448E-3   &    2.11     &  3.1546E-3   &   2.10      &  1.1120E-2   &   2.04     \\  
        %4     &  5.5512E-4   &    2.05     &  7.4180E-4   &    2.08     &  7.4438E-4   &   2.08      &  2.6914E-3   &   2.05     \\  
        &\multicolumn{8}{c}{ROD-ADER-DG-$\mathcal{P}_2$}\\[0.5mm]
        0     &  2.1023E-2   &    --       &  2.6769E-2   &    --       &  2.8488E-2   &   --        &  9.5588E-2   &   --       \\  
        1     &  5.4789E-3   &    1.94     &  7.1941E-3   &    1.90     &  7.2213E-3   &   1.98      &  2.6524E-2   &   1.85     \\  
        2     &  1.2612E-3   &    2.12     &  1.6497E-3   &    2.12     &  1.6571E-3   &   2.12      &  6.1478E-3   &   2.11     \\  
        3     &  2.9615E-4   &    2.09     &  3.9261E-4   &    2.07     &  3.9366E-4   &   2.07      &  1.4588E-3   &   2.07     \\  
        %4     &  7.0454E-5   &    2.07     &  9.4221E-5   &    2.06     &  9.4451E-5   &   2.06      &  3.4882E-4   &   2.06     \\  
        &\multicolumn{8}{c}{ROD-ADER-DG-$\mathcal{P}_3$}\\[0.5mm]
        0     &  1.2482E-3   &    --       &  1.4765E-3   &    --       &  1.4488E-3   &   --        &  4.7772E-3   &   --       \\  
        1     &  5.7882E-5   &    4.43     &  8.1071E-5   &    4.18     &  8.0220E-5   &   4.17      &  2.3357E-4   &   4.35     \\  
        2     &  6.6352E-6   &    3.12     &  6.3173E-6   &    3.68     &  6.2445E-6   &   3.68      &  2.3344E-5   &   3.32     \\  
        3     &  1.1183E-6   &    2.57     &  1.1242E-6   &    2.49     &  1.1643E-6   &   2.42      &  4.4317E-6   &   2.40     \\  
        \hline\\[1pt]
\end{tabular}
\end{table}

\subsection{Shock-cylinder interaction}\label{test:shkcyl}

\begin{figure}
\centering
\subfloat[Conformal]{\includegraphics[height=0.4\textwidth]{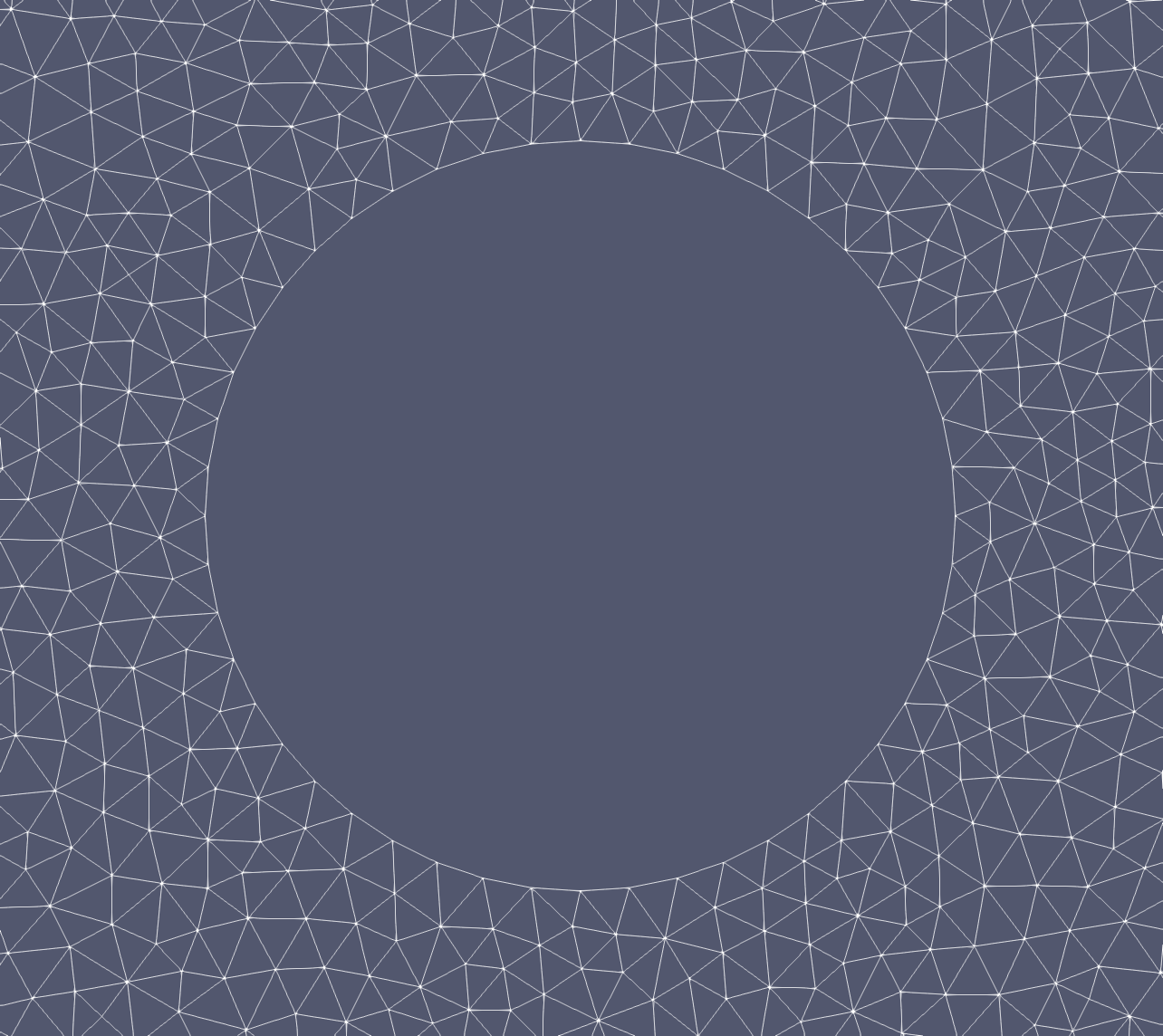}}\quad
\subfloat[Embedded]{\includegraphics[height=0.4\textwidth]{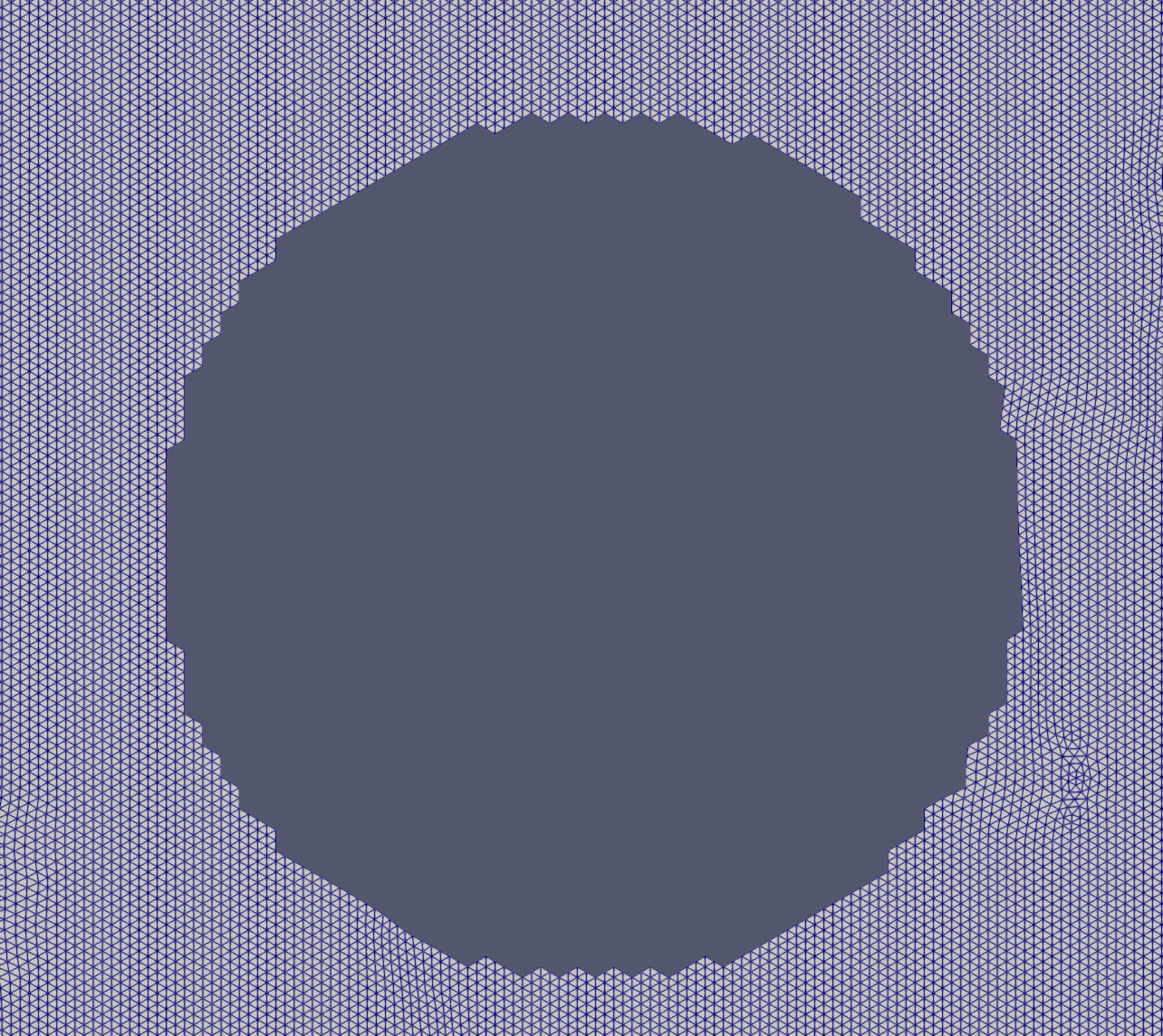}}
\caption{Shock-cylinder interaction: characteristics of computational meshes on a computational domain $[-2,6]\times[-3,3]$ for the test case presented in~\cref{test:shkcyl}.}\label{fig:shkcylmesh} 
\end{figure}

In the last test case, we present the possibility of tackling compressible flows featuring discontinuities with a 
consistent ROD boundary treatment to improve the accuracy close to internal boundaries.  
The employed ADER-DG scheme is supplemented with the \textit{a posteriori} sub-cell finite volume limiter discussed in~\cref{sec:limiter}.
Here we present the numerical results obtained on a computational domain $[-2,6]\times[-3,3]$ discretized
with a unstructured conformal triangulation and a fully embedded one (see~\cref{fig:shkcylmesh}).
The cylinder is centered in $(0,0)$ and has radius 0.5.
The initial condition consists in a shock wave traveling at Mach number $M_s = 1.3$ and is then given via the Rankine-Hugoniot conditions.
The flow upstream the shock is at rest and is characterized by density and pressure, respectively being $\rho=1.4$ and $p=1$.

\begin{figure}
\centering
\subfloat[with ROD ($t=1.0$)]{\includegraphics[height=0.35\textwidth,trim={10cm 8cm 26cm 8cm},clip]{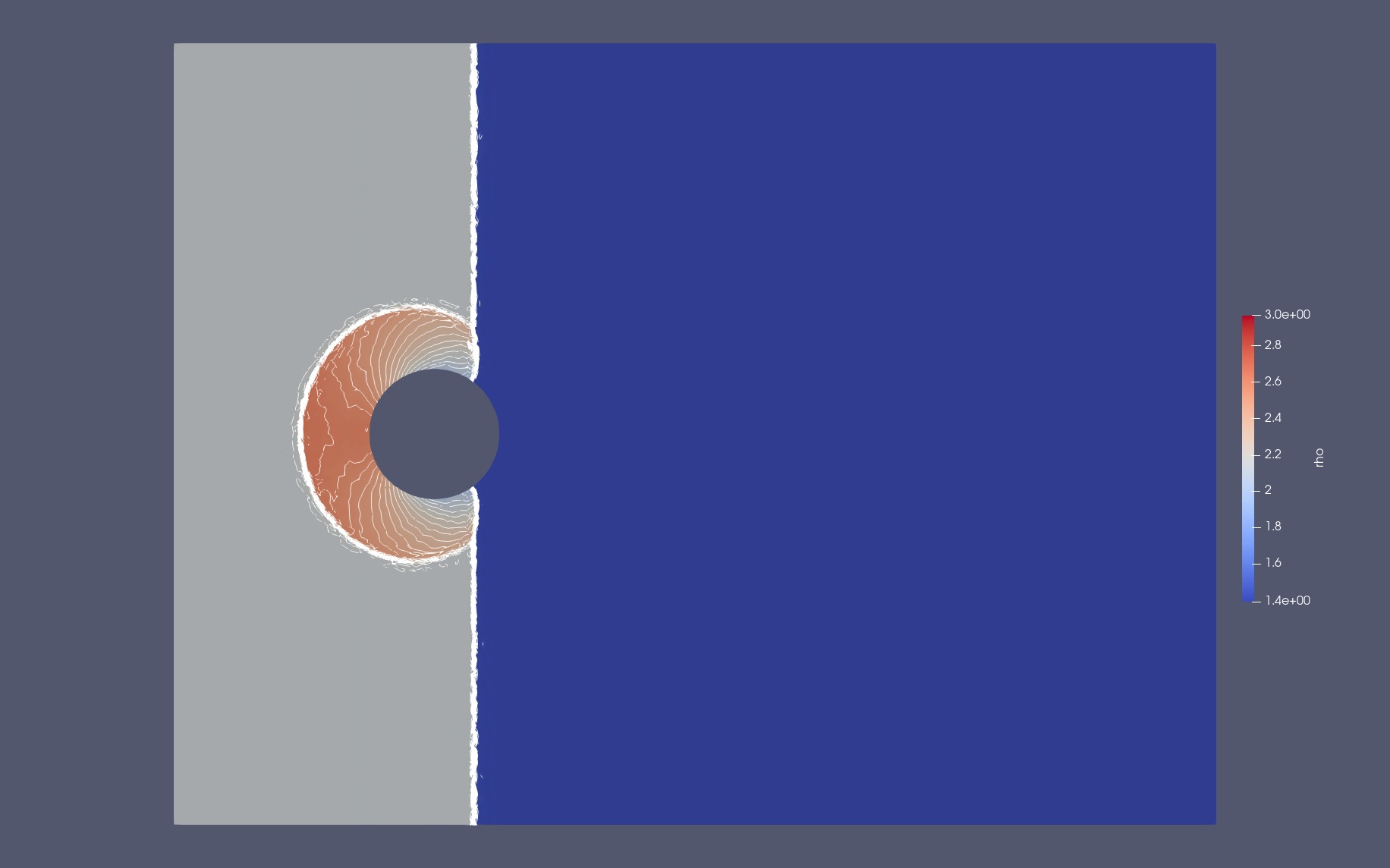}}\quad
\subfloat[without ROD ($t=1.0$)]{\includegraphics[height=0.35\textwidth,trim={10cm 8cm 26cm 8cm},clip]{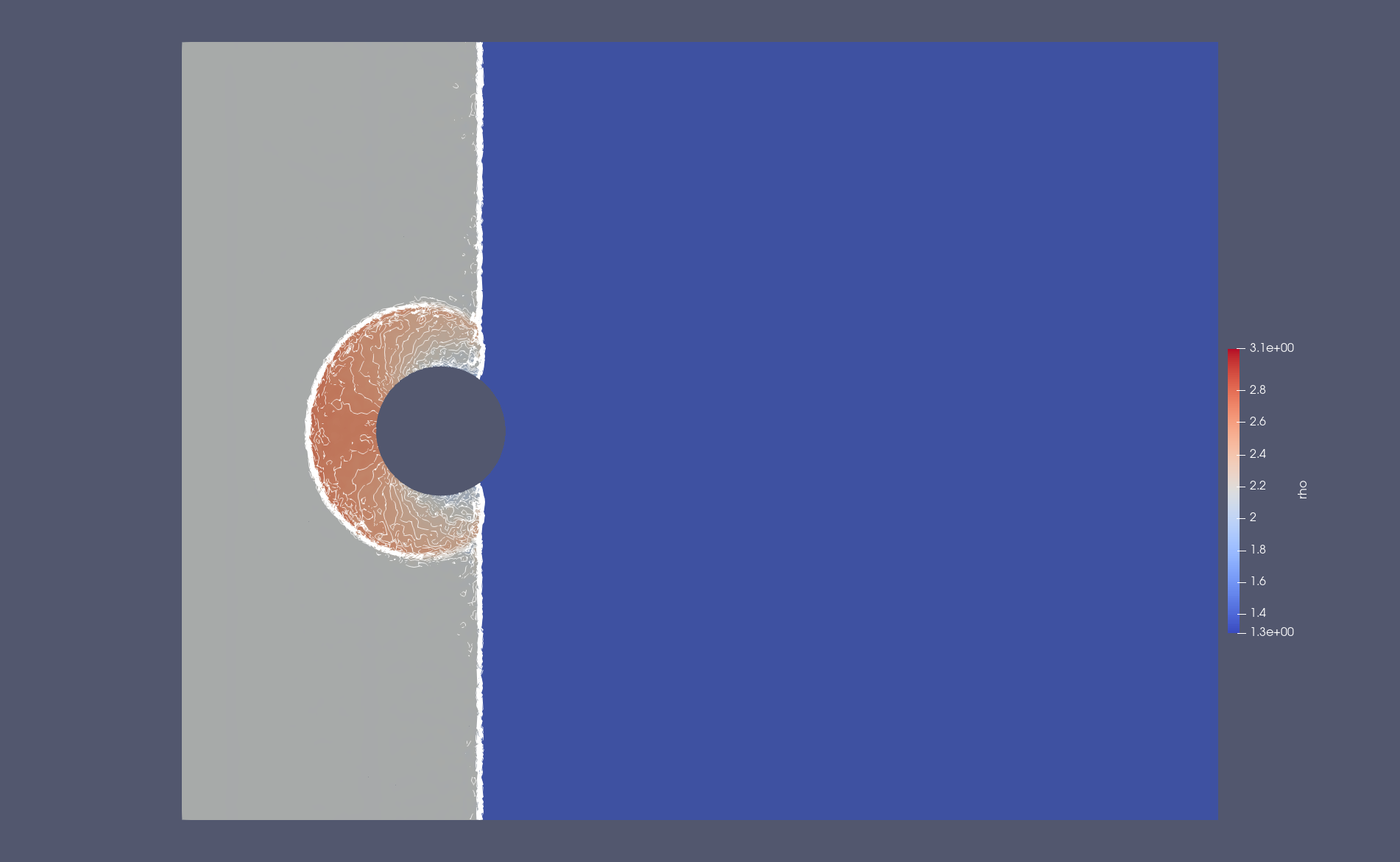}} \vspace{1em}
\subfloat[with ROD ($t=1.5$)]{\includegraphics[height=0.35\textwidth,trim={10cm 8cm 26cm 8cm},clip]{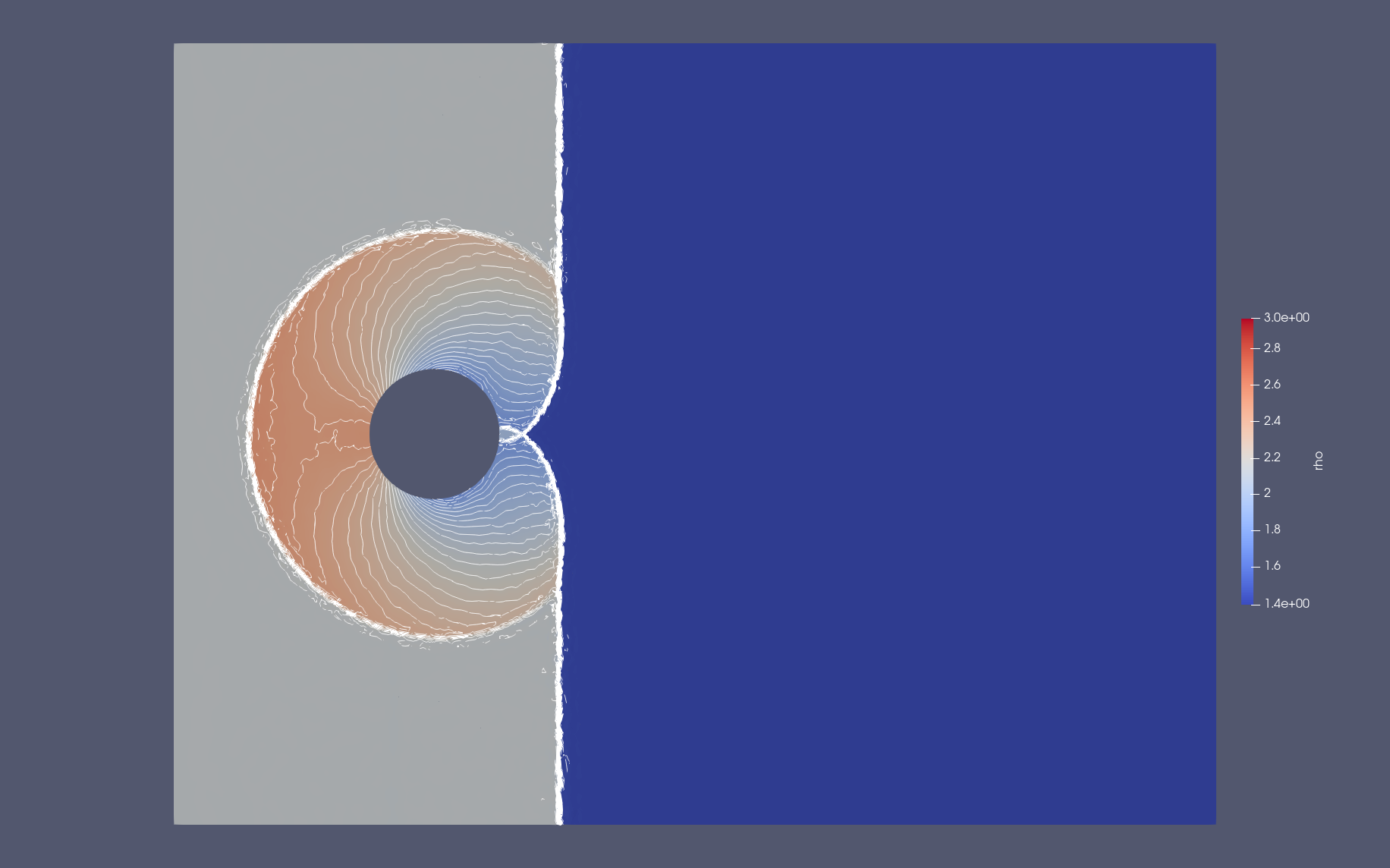}}\quad
\subfloat[without ROD ($t=1.5$)]{\includegraphics[height=0.35\textwidth,trim={10cm 8cm 26cm 8cm},clip]{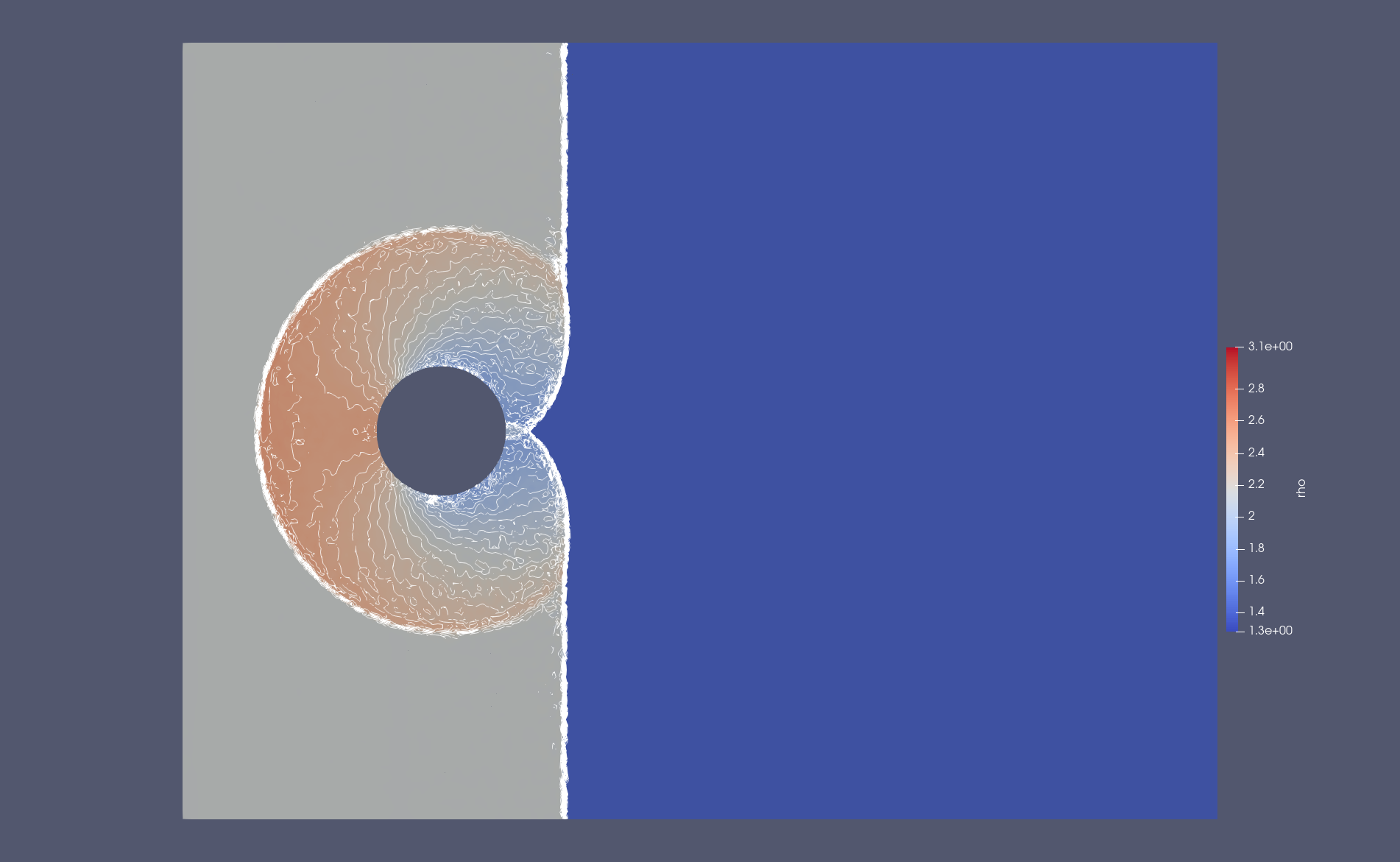}} \vspace{1em}
\subfloat[with ROD ($t=2.0$)]{\includegraphics[height=0.35\textwidth,trim={9.0cm 8cm 27cm 8cm},clip]{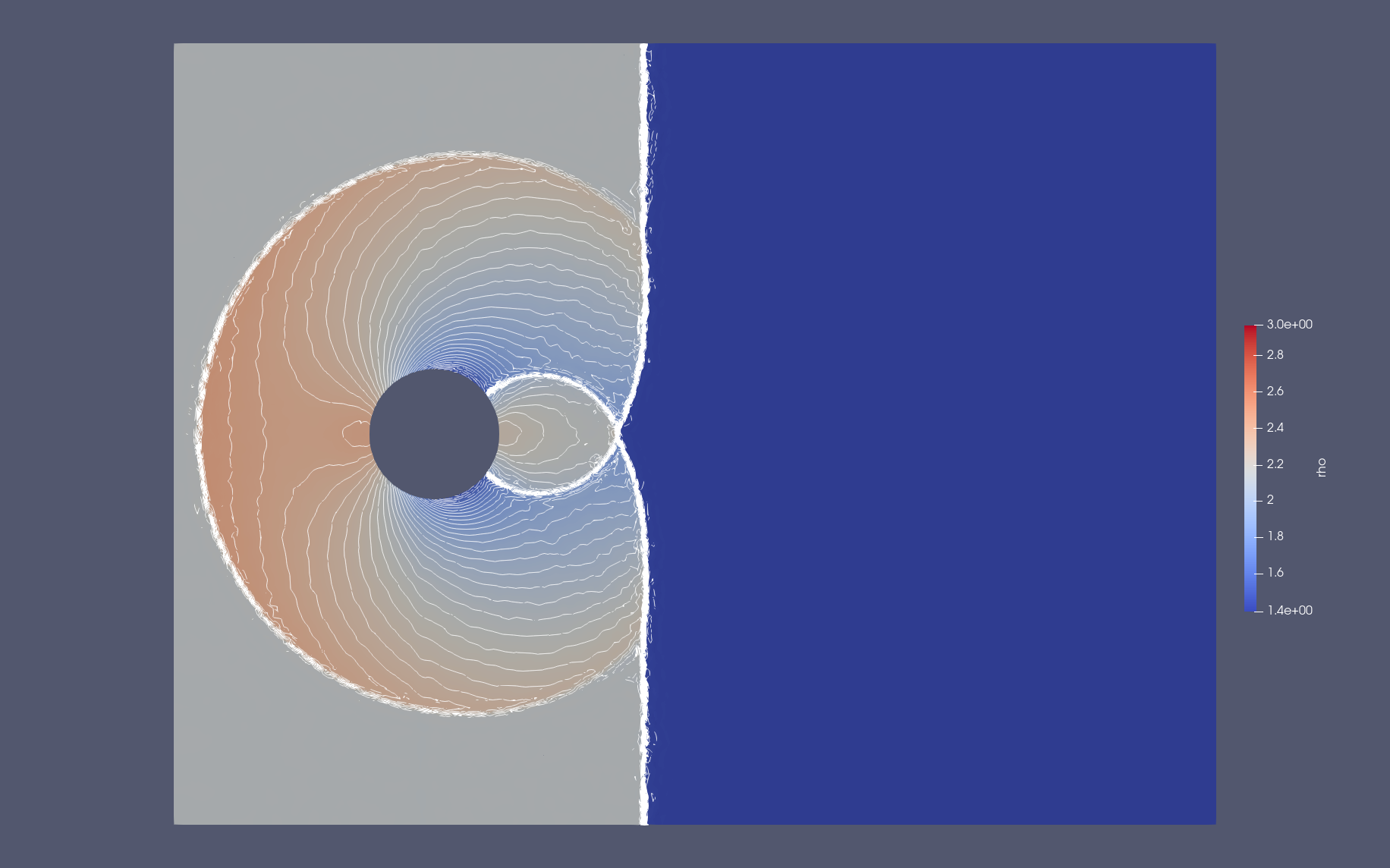}}\quad
\subfloat[without ROD ($t=2.0$)]{\includegraphics[height=0.35\textwidth,trim={10cm 8cm 26cm 8cm},clip]{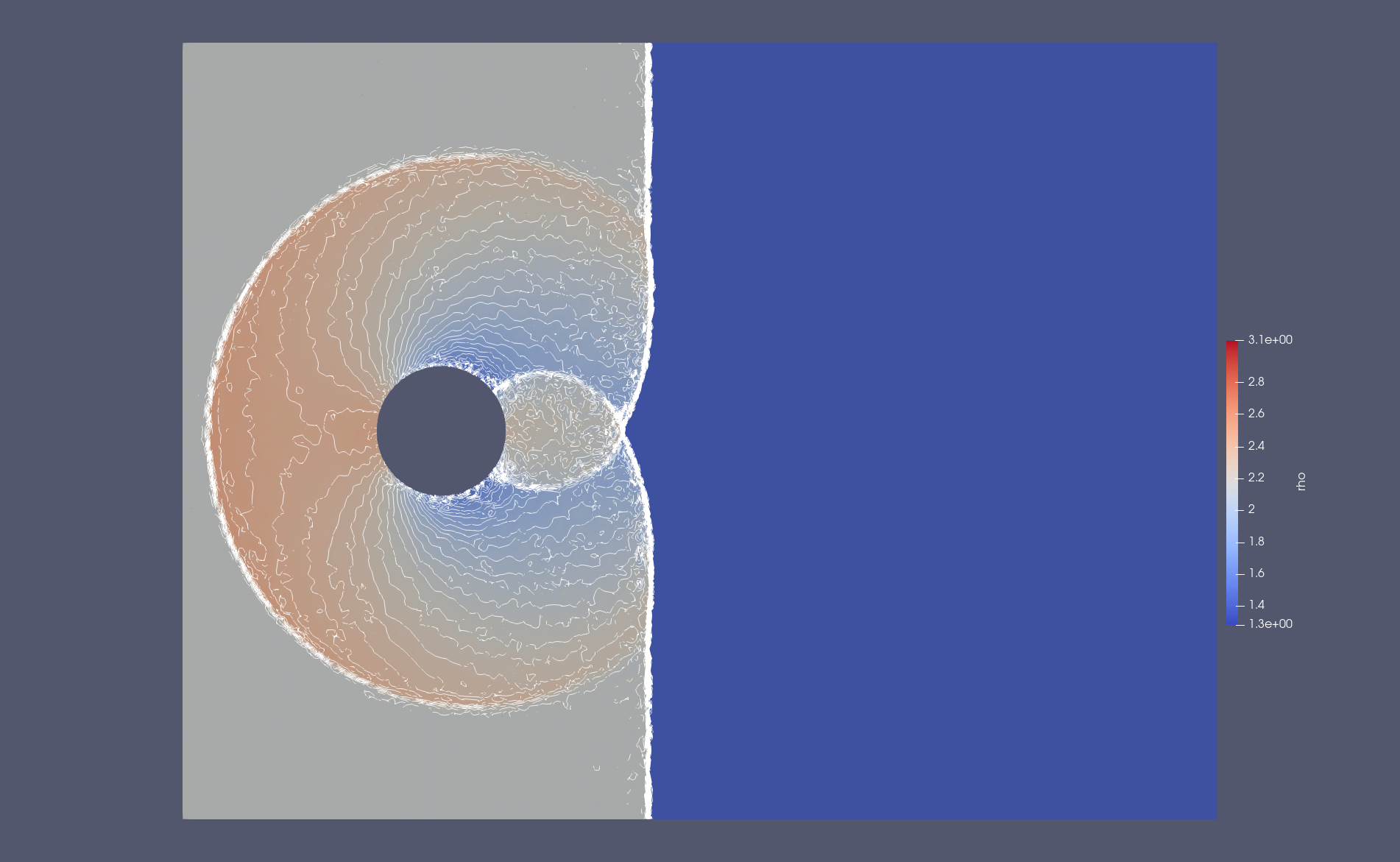}}
\caption{Shock-cylinder interaction: density iso-contours at different times for the test case presented in~\cref{test:shkcyl}. We show the results obtained with classical reflecting boundary conditions (right), and the wall conditions enforced by the ROD boundary polynomial (left) on the conformal mesh presented in~\cref{fig:shkcylmesh}a.}\label{fig:shkcylconformalcomparison}
\end{figure}

In~\cref{fig:shkcylconformalcomparison} we present the numerical results obtained with third order polynomials 
on the conformal triangulation by imposing naive reflecting boundary conditions (without ROD) and the 
general slip conditions (with ROD). It is clear that when curved boundaries are considered, a more general boundary treatment 
that takes into account the real geometry is essential to capture properly the flow dynamics.  
In the simulation performed without ROD, we notice that spurious disturbances starts from the poor boundary treatment and propagates
in the computational domain making the solution less accurate. 
On the contrary, much smoother iso-contours both close and away from the boundary are observed in the simulation performed with the ROD
method.

\begin{figure}
\centering
\subfloat[Density ($t=1.0$)]{\includegraphics[height=0.35\textwidth]{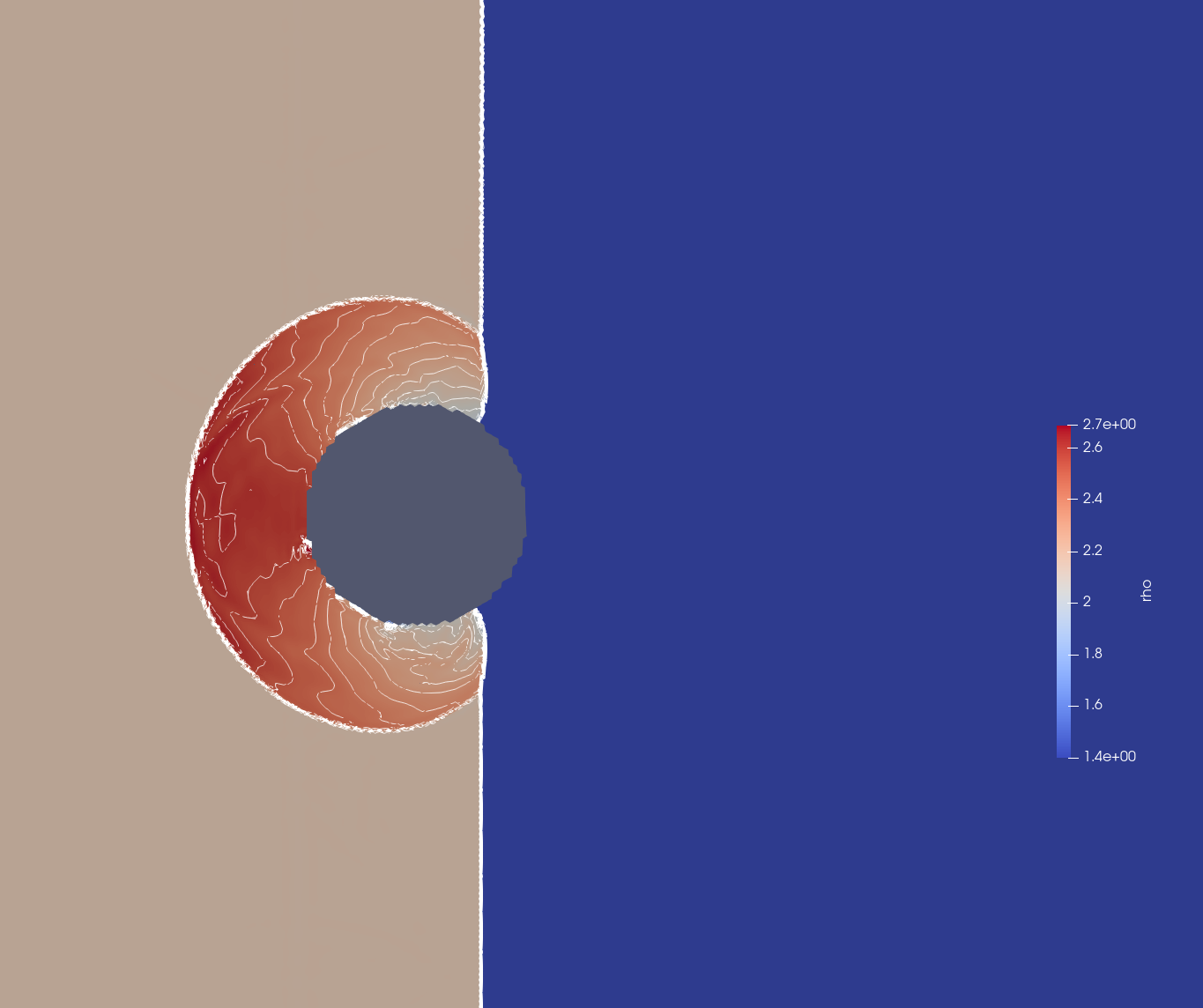}}\quad
\subfloat[Limiter ($t=1.0$)]{\includegraphics[height=0.35\textwidth]{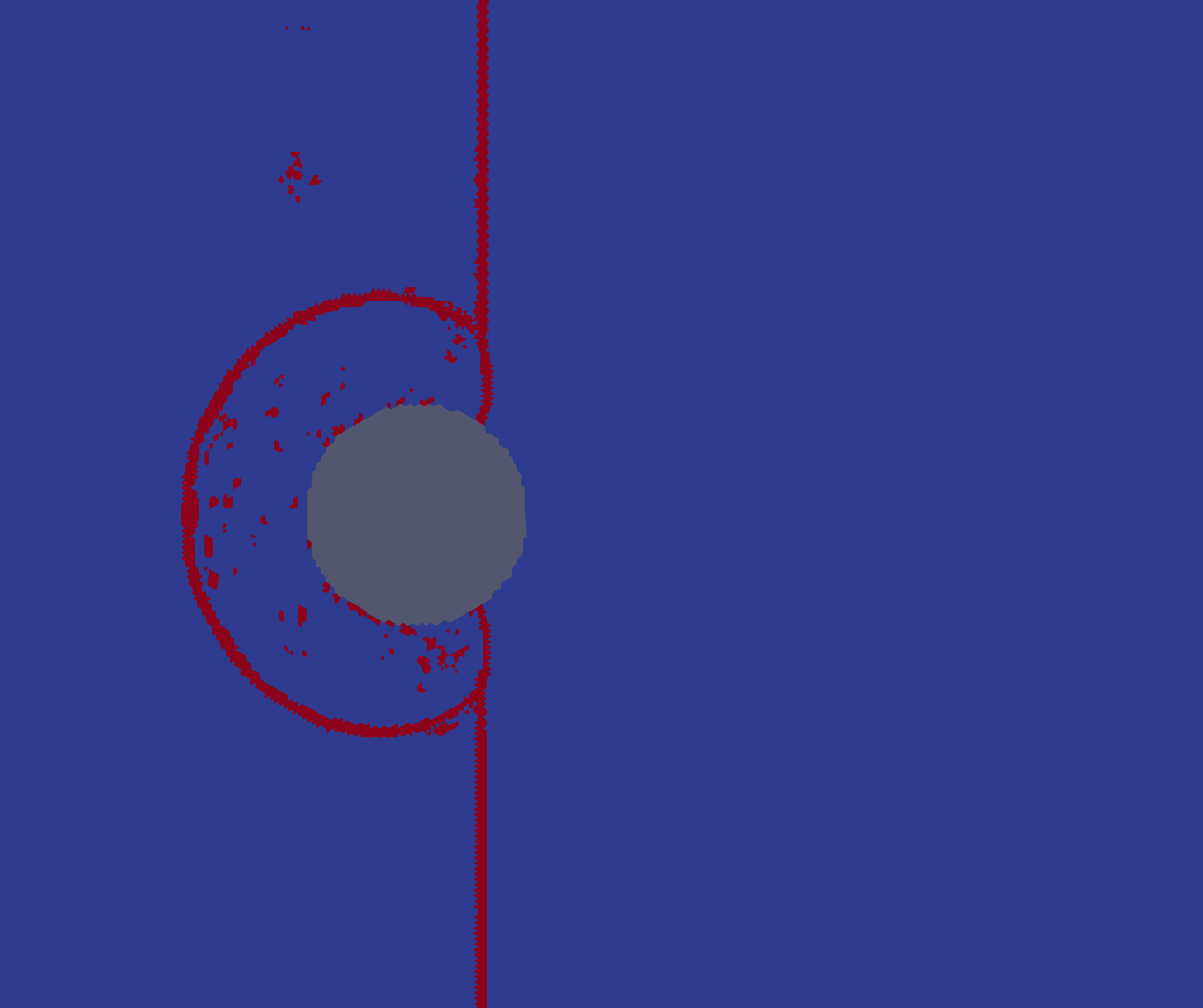}}\qquad
\subfloat[Density ($t=1.5$)]{\includegraphics[height=0.35\textwidth]{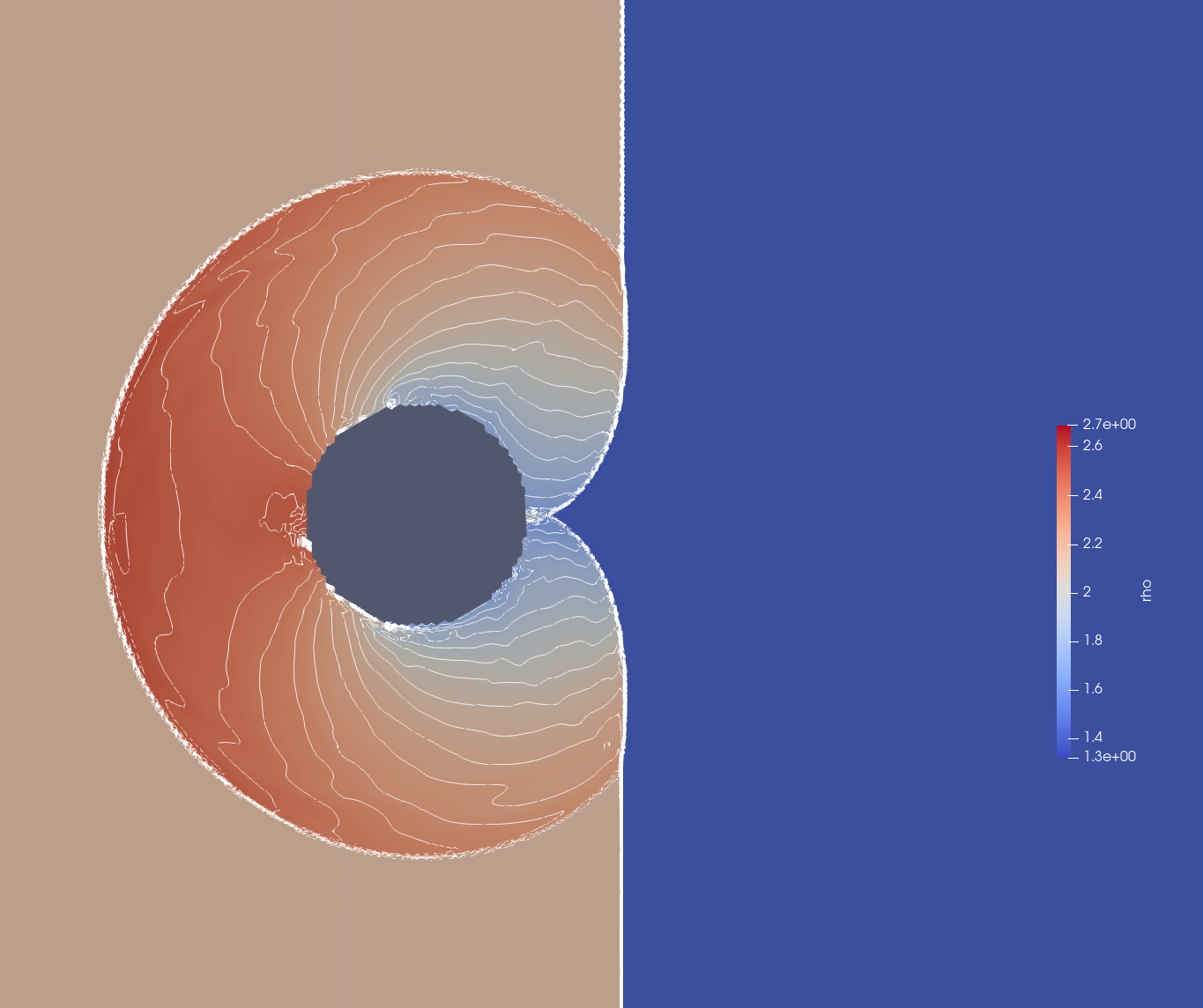}}\quad
\subfloat[Limiter ($t=1.5$)]{\includegraphics[height=0.35\textwidth]{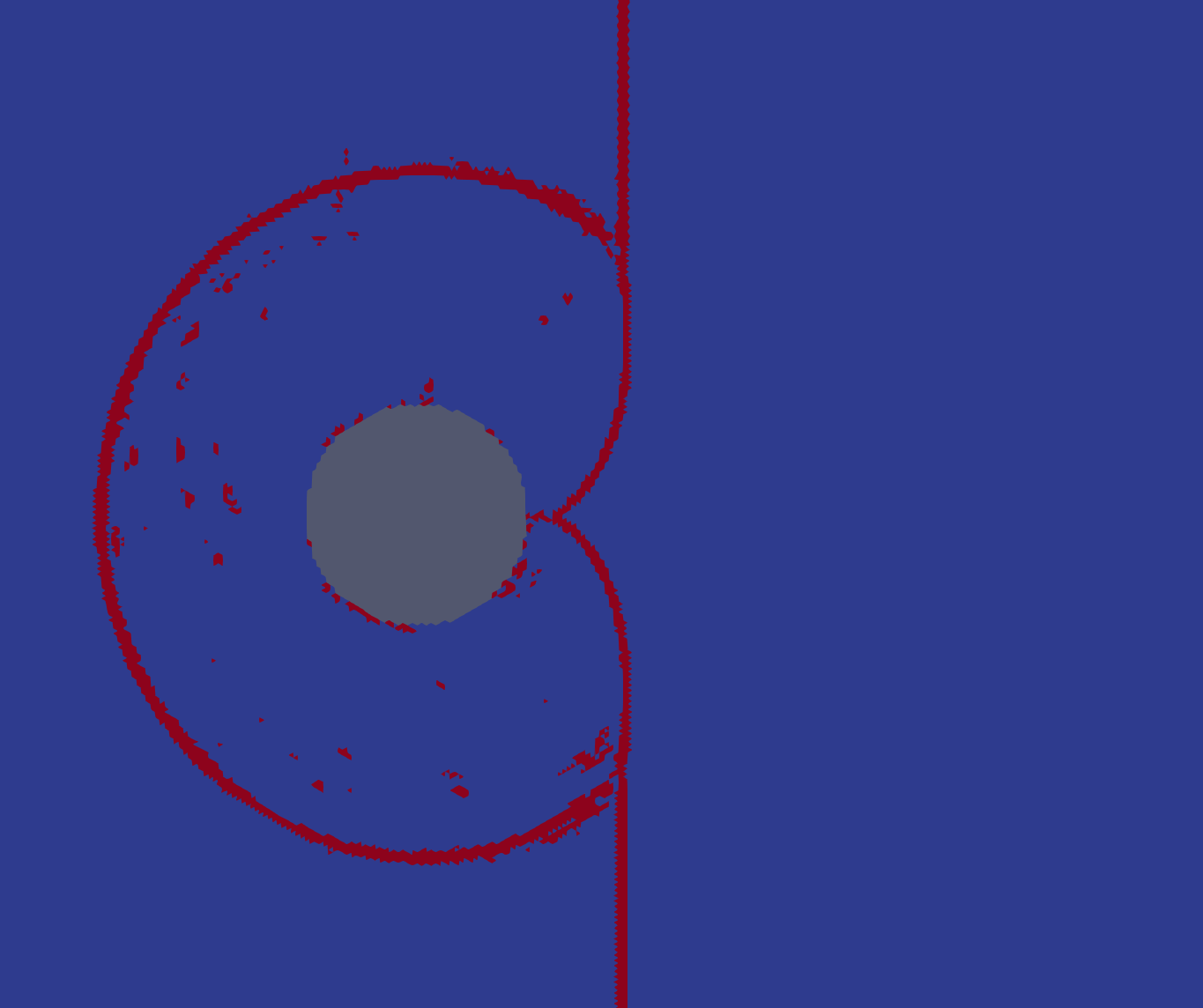}}\qquad
\subfloat[Density ($t=2.0$)]{\includegraphics[height=0.35\textwidth]{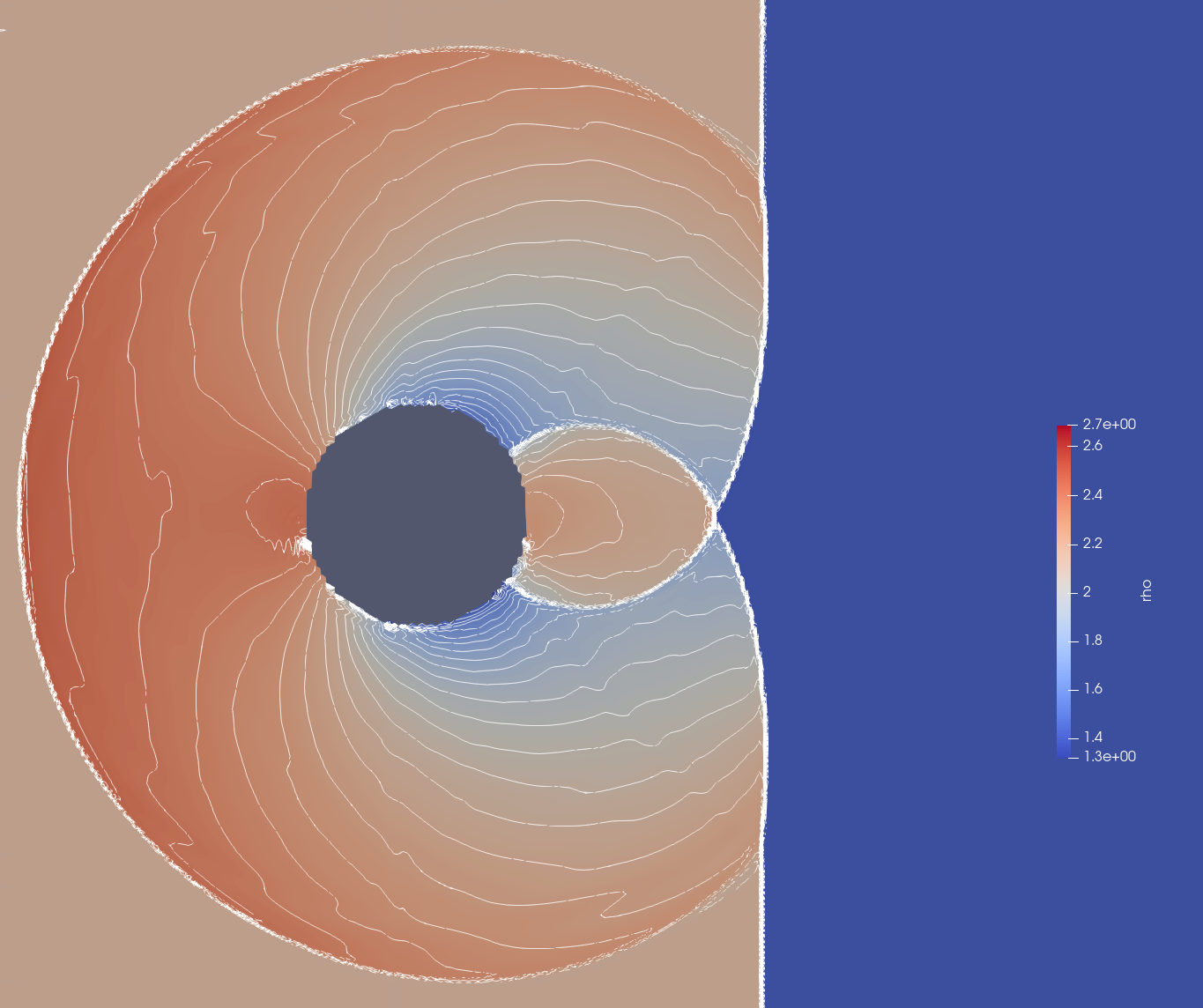}}\quad
\subfloat[Limiter ($t=2.0$)]{\includegraphics[height=0.35\textwidth]{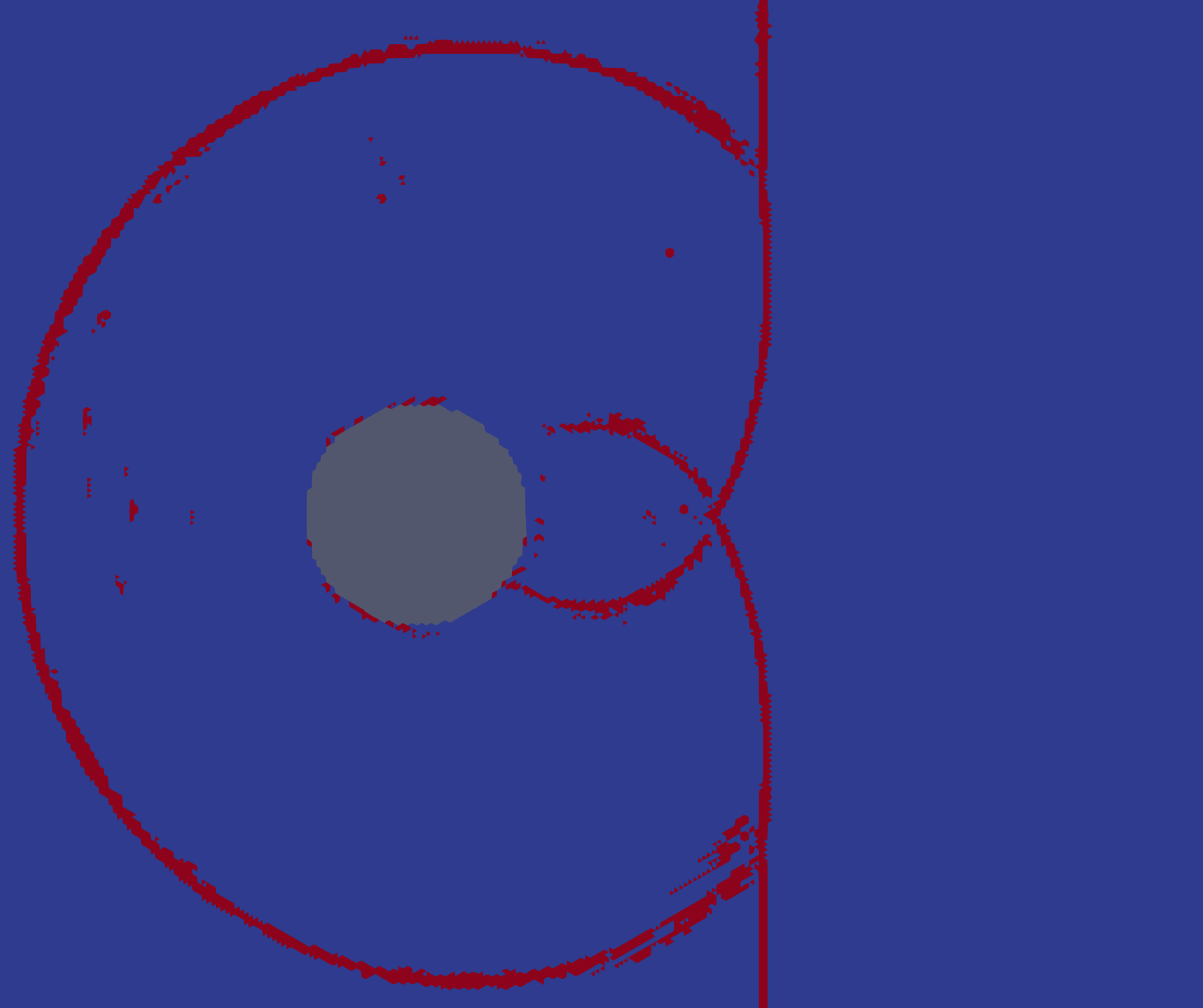}}     
\caption{Shock-cylinder interaction: density iso-contours and troubled cells at different times for the test case presented in~\cref{test:shkcyl}. We show the results obtained with the ROD boundary polynomial on the fully embedded mesh presented in~\cref{fig:shkcylmesh}b.}\label{fig:shkcylembedded}
\end{figure}

Finally, in~\cref{fig:shkcylembedded} a numerical simulation run on a fully embedded triangulation is 
presented. Here we show the results obtained by imposing the general slip condition by means of the ROD boundary reconstruction on
a surrogate boundary completely detached from the exact geometry as shown in~\cref{fig:shkcylmesh}. 
In general the density iso-contours appears globally very similar to the one obtained in the conformal case performed with the ROD method.
Some spurious disturbances present in some cells of the discretized boundary are normal due to the use of a less accurate 
boundary treatment when the {\it a posteriori} limiter is activated. For this reason, we present side by side also the troubled cells
recomputed by means of the low order method to show that disturbances only arise in those regions, while in the non-troubled parts
where the ROD treatment is applied the solution appears smooth.

\section{Conclusions and future investigations}
\label{sec:conclusions}

In this work, we have investigated the numerical aspects of the Reconstruction for Off-site Data (ROD) method in a 
discontinuous Galerkin and space-time ADER framework. 
The ROD method is formulated as a constraint minimization problem that allows to recover consistent boundary polynomials 
used for the correct imposition of boundary conditions. Originally developed in a finite volume framework, the ROD idea consisted 
in managing multiple constraints related to scattered mean values associated with neighboring elements, found with specific algorithms.
Thanks to the local discontinuous finite element formulation, we show how the ROD method can be easily adapted by evaluating the 
boundary cell polynomial within quadrature points, rather than using neighboring elements. By doing so, a classical constrained least squares
should be solved with Lagrange multipliers for the boundary integral of the semi-discretization~\eqref{eq:DGweak semidiscrete 2}.  
Naturally, it follows that by integrating~\eqref{eq:DGweak semidiscrete 2} with classical high order time integration schemes, 
like one-step Runge-Kutta methods, such minimization problems must be solved at each sub-time node.
Although the linear system arising from the constrained least squares is local, an inversion is needed at each boundary cell for 
each sub-time step, hence impacting the computational costs. 
For instance, when dealing with $\mathcal P_2$ polynomials a local matrix of $9\times 9$ must be inverted each time;
for $\mathcal P_3$ the matrix size becomes $14\times 14$. 
In this paper, we also explored a new advantageous perspective for the ROD method when considering the space-time ADER framework.
Indeed, we showed that by using the same space-time polynomial coefficients of the ADER predictor~\eqref{eq:predictor0} 
for the minimization problem, we obtain a simplified linear system that provides a  polynomial 
consistent with the physical boundary condition for the entire space-time element, hence only one inversion is needed. 
Several numerical experiments are shown to assess the high order convergence properties of both the spatial and space-time ROD
reconstructions when dealing with Dirichlet conditions on a curved boundary using conformal linear meshes and fully embedded triangulations. 
Application of the novel space-time formulation to general Dirichlet-Neumann-Robin boundary conditions is also presented 
for the general slip-wall treatment. The latter opens up to many applications when the interaction between a flow and a boundary is crucial.
In the last test case, the application of the ROD method to a shock-cylinder interaction is presented, where the discontinuous features of
the flow are properly captured using an {\it a posteriori} MOOD approach.

Several important aspects related to the proposed numerical method need further investigations. 
In particular, we believe that the introduction of additional constraints within the optimization problem may bring
more accurate results and improve the conditioning of the problem and the activation of the shock limiter on boundary cells. 
For instance, mass conservation, positivity and entropy preservation will be considered in future works.

\section*{Acknowledgments}
M.R.\ is a member of the CARDAMOM team, INRIA and University of Bordeaux research center.
E.G.\ gratefully acknowledges the support and funding received from the European Union with the ERC Starting Grant ALcHyMiA (No. 101114995).  
Views and opinions expressed are however those of the author(s) only and do not necessarily reflect those of the European Union or the European Research Council Executive Agency. 
Neither the European Union nor the granting authority can be held responsible for them.

\bibliographystyle{elsarticle-num} 
\bibliography{biblio.bib}

\end{document}